\documentclass[11pt]{article}
\usepackage{amsfonts}
\usepackage{amsmath}
\usepackage[a4paper]{geometry}
\usepackage[onehalfspacing]{setspace}
\usepackage{graphicx}

\newtheorem{theorem}{Theorem}

\newtheorem{corollary}[theorem]{Corollary}

\newenvironment{proof}[1][Proof]{\noindent\textbf{#1.} }{\ \rule{0.5em}{0.5em}}
\input{tcilatex}
\begin{document}

\title{Generalised vector products in three-dimensional geometry}
\author{G A Notowidigdo* \\
School of Mathematics and Statistics\\
UNSW Sydney\\
Sydney,\ NSW, Australia\\
gnotowidigdo@zohomail.com.au \and N J Wildberger \\
School of Mathematics and Statistics\\
UNSW Sydney\\
Sydney, NSW, Australia\\
n.wildberger@unsw.edu.au}
\date{}
\maketitle

\begin{abstract}
In three-dimensional Euclidean geometry, the scalar product produces a
number associated to two vectors, while the vector product computes a vector
perpendicular to them. These are key tools of physics, chemistry and
engineering and supported by a rich vector calculus of 18th and 19th century
results. This paper extends this calculus to arbitrary metrical geometries
on three-dimensional space, generalising key results of Lagrange, Jacobi,
Binet and Cauchy in a purely algebraic setting which applies also to general
fields, including finite fields. We will then apply these vector theorems to
set up the basic framework of rational trigonometry in the three-dimensional
affine space and the related two-dimensional projective plane, and show an
example of its applications to relativistic geometry.
\end{abstract}



\section{Introduction}

In three-dimensional Euclidean vector geometry, the \textbf{scalar product}
of two vectors 
\begin{equation*}
v_{1}\equiv \left( x_{1},y_{1},z_{1}\right) \quad \mathrm{and}\quad
v_{2}\equiv \left( x_{2},y_{2},z_{2}\right)
\end{equation*}%
is the number%
\begin{equation*}
v_{1}\cdot v_{2}\equiv x_{1}x_{2}+y_{1}y_{2}+z_{1}z_{2}
\end{equation*}%
while the \textbf{vector product} is the vector%
\begin{equation*}
v_{1}\times v_{2}=\left( x_{1},y_{1},z_{1}\right) \times \left(
x_{2},y_{2},z_{2}\right) \equiv \left(
y_{1}z_{2}-y_{2}z_{1},x_{2}z_{1}-x_{1}z_{2},x_{1}y_{2}-x_{2}y_{1}\right) .
\end{equation*}%
These definitions date back to the late eighteenth century, with Lagrange's 
\cite{Lagrange} study of the tetrahedron in three-dimensional space.
Hamilton's quaternions in four-dimensional space combine both operations,
but Gibbs \cite{Gibbs} and Heaviside \cite{Heaviside} choose to separate the
two quantities and introduce the current notations, where they were
particularly effective in simplifying Maxwell's equations.

Implicitly the framework of Euclidean geometry from Euclid's \textit{Elements%
} \cite{Heath} underlies both concepts, and gradually in the 20th century it
became apparent that the scalar, or dot, product in particular could be
viewed as the starting point for Euclidean metrical geometry, at least in a
vector space. However the vector, or cross, product is a special
construction more closely aligned with three-dimensional space, although the
framework of Geometric Algebra (see \cite{Hestenes et al} and \cite{Doran et
al}) does provide a more sophisticated framework for generalizing it.

In this paper, we initiate the study of vector products in general
three-dimensional inner product spaces over arbitrary fields, not of
characteristic $2$, and then apply this to develop the main laws of rational
trigonometry for triangles both in the affine and projective settings, now
with respect to a general symmetric bilinear form, and over an arbitrary
field.

We will start with the associated three-dimensional vector space $\mathbb{V}%
^{3}$, regarded as row vectors, of a three-dimensional affine space $\mathbb{%
A}^{3}$ over an arbitrary field $\mathbb{F}$, not of characteristic $2.$ A
general metrical framework is introduced via a $3\times 3$ invertible
symmetric matrix $B$, which defines a non-degenerate \textbf{symmetric
bilinear form} on the vector space $\mathbb{V}^{3}$ by%
\begin{equation*}
v\cdot _{B}w\equiv vBw^{T}.
\end{equation*}%
We will call this number, which is necessarily an element of the field $%
\mathbb{F},$ the $B$\textbf{-scalar product} of $v$ and $w$.

If $B$ is explicitly the symmetric matrix 
\begin{equation*}
B\equiv 
\begin{pmatrix}
a_{1} & b_{3} & b_{2} \\ 
b_{3} & a_{2} & b_{1} \\ 
b_{2} & b_{1} & a_{3}%
\end{pmatrix}%
\end{equation*}%
then the \textbf{adjugate }of\textbf{\ }$B$ is defined by 
\begin{equation*}
\limfunc{adj}B=%
\begin{pmatrix}
a_{2}a_{3}-b_{1}^{2} & b_{1}b_{2}-a_{3}b_{3} & b_{1}b_{3}-a_{2}b_{2} \\ 
b_{1}b_{2}-a_{3}b_{3} & a_{1}a_{3}-b_{2}^{2} & b_{2}b_{3}-a_{1}b_{1} \\ 
b_{1}b_{3}-a_{2}b_{2} & b_{2}b_{3}-a_{1}b_{1} & a_{1}a_{2}-b_{3}^{2}%
\end{pmatrix}%
\end{equation*}%
and in case $B$ is invertible, 
\begin{equation*}
\limfunc{adj}B=\left( \det B\right) B^{-1}.
\end{equation*}%
Then we define the $B$\textbf{-vector product} between the two vectors $%
v_{1} $ and $v_{2}$ to be the vector 
\begin{equation*}
v_{1}\times _{B}v_{2}\equiv \left( v_{1}\times v_{2}\right) \func{adj}B.
\end{equation*}

We will extend well-known and celebrated formulas from Binet \cite{Binet},
Cauchy \cite{Cauchy}, Jacobi \cite{Jacobi} and Lagrange \cite{Lagrange}
involving also the $B$\textbf{-scalar triple product }%
\begin{equation*}
\left[ v_{1},v_{2},v_{3}\right] _{B}\equiv v_{1}\cdot _{B}\left( v_{2}\times
_{B}v_{3}\right)
\end{equation*}%
the $B$-\textbf{vector triple product} 
\begin{equation*}
\left\langle v_{1},v_{2},v_{3}\right\rangle _{B}\equiv v_{1}\times
_{B}\left( v_{2}\times _{B}v_{3}\right)
\end{equation*}%
the $B$\textbf{-scalar quadruple product} 
\begin{equation*}
\left[ v_{1},v_{2};v_{3},v_{4}\right] _{B}\equiv \left( v_{1}\times
_{B}v_{2}\right) \cdot _{B}\left( v_{3}\times _{B}v_{4}\right)
\end{equation*}%
and the $B$-\textbf{vector quadruple product} 
\begin{equation*}
\left\langle v_{1},v_{2};v_{3},v_{4}\right\rangle _{B}\equiv \left(
v_{1}\times _{B}v_{2}\right) \times _{B}\left( v_{3}\times _{B}v_{4}\right) .
\end{equation*}

These tools will then be the basis for a rigorous framework of vector
trigonometry over the three-dimensional vector space $\mathbb{V}^{3}$,
following closely the framework of affine rational trigonometry in
two-dimensional Euclidean space formulated in \cite{WildDP}, but working
more generally with \textit{vector triangles }in three-dimensional vector
space, by which we mean an unordered collection of three vectors whose sum
is the zero vector. In this framework, the terms "quadrance" and "spread"
are used rather than "distance" and "angle" respectively; though in
Euclidean geometry quadrances and spreads are equivalent to squared
distances and squared sines of angles respectively, when we start to
generalise to arbitrary symmetric bilinear forms on vector spaces over
general fields the definitions of quadrance and spread easily generalise
using only linear algebraic methods. To highlight this in the case of
Euclidean geometry over a general field, the quadrance of a vector $v$ is
the number%
\begin{equation*}
Q\left( v\right) \equiv v\cdot v
\end{equation*}%
and the spread between two vectors $v$ and $w$ is the number%
\begin{equation*}
s\left( v,w\right) \equiv 1-\frac{\left( v\cdot w\right) ^{2}}{\left( v\cdot
v\right) \left( w\cdot w\right) }.
\end{equation*}%
By using these quantities, we avoid the issue of taking square roots of
numbers which are not square numbers in such a general field, which becomes
problematic when working over finite fields.

We also develop a corresponding framework of \textit{projective triangles }%
in the two-dimensional projective space associated to the three-dimensional
vector space, where geometrically a projective triangle can be viewed as a
triple of one-dimensional subspaces of $\mathbb{V}^{3}$. Mirroring the setup
in \cite{WildAP} and \cite{WildUHG1}, the tools presented above will be
sufficient to develop the metrical geometry of projective triangles and
extend known results from \cite{WildAP} and \cite{WildUHG1} to arbitrary
non-degenerate $B$-scalar products. This then gives us a basis for both
general elliptic and hyperbolic rational geometries, again working over a
general field.

Applying the tools of $B$-scalar and $B$-vector products together gives a
powerful method in which to study two-dimensional rational trigonometry in
three-dimensional space over a general metrical framework and ultimately the
tools presented in the paper could be used to study three-dimensional
rational trigonometry in three-dimensional space over said general metrical
framework, with emphasis on the trigonometry of a general tetrahedron. This
is presented in more detail in the follow-up paper \textit{Generalised
vector products applied to affine rational trigonometry of a general
tetrahedron}.

We will illustrate this new technology with several explicit examples where
all the calculations are completely visible and accurate, and can be done by
hand. First we look at the geometry of a methane molecule $CH_{4}$ which is
a regular tetrahedron, and derive new (rational!) expressions for the
separation of the faces, and the solid spreads. We also analyse both a
vector triangle and a projective triangle in a vector space with a Minkowski
bilinear form; this highlights how the general metrical framework of affine
and projective rational trigonometry can be used to study relativistic
geometry.

\section{Vector algebra over a general metrical framework}

We start by considering the three-dimensional vector space $\mathbb{V}^{3}$
over a general field $\mathbb{F}$ not of characteristic $2$, consisting of
row vectors 
\begin{equation*}
v=\left( x,y,z\right) .
\end{equation*}

\subsection{The $B$-scalar product}

A $3\times 3$ symmetric matrix%
\begin{equation}
B\equiv 
\begin{pmatrix}
a_{1} & b_{3} & b_{2} \\ 
b_{3} & a_{2} & b_{1} \\ 
b_{2} & b_{1} & a_{3}%
\end{pmatrix}
\label{Symmetric Matrix}
\end{equation}%
determines a \textbf{symmetric bilinear form}\ on $\mathbb{V}^{3}$ defined by%
\begin{equation*}
v\cdot _{B}w\equiv vBw^{T}\mathbf{.}
\end{equation*}%
We will call this the $B$-\textbf{scalar product}. The $B$-scalar product is 
\textbf{non-degenerate} precisely when the condition that $v\cdot _{B}w=0$
for any vector $v$ in $\mathbb{V}^{3}$ implies that $w=0$; this will occur
precisely when $B$ is invertible, that is when%
\begin{equation*}
\det B\neq 0.
\end{equation*}%
We will assume that the $B$-scalar product is non-degenerate throughout this
paper.

The associated $B$-\textbf{quadratic form }on $\mathbb{V}^{3}$ is defined by%
\textbf{\ }%
\begin{equation*}
Q_{B}\left( v\right) \equiv v\cdot _{B}v
\end{equation*}%
and the number $Q_{B}\left( v\right) $ is the $B$\textbf{-quadrance }of%
\textbf{\ }$v$. A vector $v$ is $B$\textbf{-null }precisely when 
\begin{equation*}
Q_{B}\left( v\right) =0.
\end{equation*}%
The $B$-quadrance satisfies the obvious properties that for vectors $v$ and $%
w$ in $\mathbb{V}^{3},$ and a number $\lambda $ in $\mathbb{F}$%
\begin{equation*}
Q_{B}\left( \lambda v\right) =\lambda ^{2}Q_{B}\left( v\right)
\end{equation*}%
as well as%
\begin{equation*}
Q_{B}\left( v+w\right) =Q_{B}\left( v\right) +Q_{B}\left( w\right) +2\left(
v\cdot _{B}w\right)
\end{equation*}%
and 
\begin{equation*}
Q_{B}\left( v-w\right) =Q_{B}\left( v\right) +Q_{B}\left( w\right) -2\left(
v\cdot _{B}w\right) .
\end{equation*}%
Hence the $B$-scalar product can be expressed in terms of the $B$-quadratic
form by either of the two \textbf{polarisation formulas} 
\begin{equation*}
v\cdot _{B}w=\frac{Q_{B}\left( v+w\right) -Q_{B}\left( v\right) -Q_{B}\left(
w\right) }{2}=\frac{Q_{B}\left( v\right) +Q_{B}\left( w\right) -Q_{B}\left(
v-w\right) }{2}.
\end{equation*}

Two vectors $v$ and $w$ in $\mathbb{V}^{3}$ are $B$-\textbf{perpendicular}
precisely when $v\cdot _{B}w=0$, in which case we write $v\perp _{B}w$.

When $B$ is the $3\times 3$ identity matrix, we have the familiar \textbf{%
Euclidean scalar product}, which we will write simply as $v\cdot w$, and the
corresponding \textbf{Euclidean quadratic form }will be written just as $%
Q\left( v\right) =v\cdot v$, which will be just the \textbf{quadrance} of $v$%
. Our aim is to extend the familiar theory of scalar and quadratic forms
from the Euclidean to the general case, and then apply these in a novel way
to establish a purely algebraic, or rational, trigonometry in three
dimensions.

\subsection{The $B$-vector product}

Given two vectors $v_{1}\equiv \left( x_{1},y_{1},z_{1}\right) $ and $%
v_{2}\equiv \left( x_{2},y_{2},z_{2}\right) $ in $\mathbb{V}^{3}$, the 
\textbf{Euclidean vector product} of $v_{1}$ and $v_{2}$ is the vector%
\begin{eqnarray*}
v_{1}\times v_{2} &=&\left( x_{2},y_{2},z_{2}\right) \times \left(
x_{2},y_{2},z_{2}\right) \\
&\equiv &\left(
y_{1}z_{2}-y_{2}z_{1},x_{2}z_{1}-x_{1}z_{2},x_{1}y_{2}-x_{2}y_{1}\right) .
\end{eqnarray*}%
We now extend this notion to the case of a general symmetric bilinear form.
Let $v_{1}$, $v_{2}$ and $v_{3}$ be vectors in $\mathbb{V}^{3}$, and let $%
M\equiv 
\begin{pmatrix}
v_{1} \\ 
v_{2} \\ 
v_{3}%
\end{pmatrix}%
$ be the matrix with these vectors as rows. We define the \textbf{adjugate}
of $M$ to be the matrix%
\begin{equation*}
\limfunc{adj}M\equiv 
\begin{pmatrix}
v_{2}\times v_{3} \\ 
v_{3}\times v_{1} \\ 
v_{1}\times v_{2}%
\end{pmatrix}%
^{T}.
\end{equation*}%
If the $3\times 3$ matrix $M$ is invertible, then the adjugate is
characterized by the equation%
\begin{equation*}
\frac{1}{\left( \det M\right) }\limfunc{adj}M\equiv M^{-1}.
\end{equation*}%
In this case the properties%
\begin{equation*}
\limfunc{adj}\left( MN\right) =\left( \limfunc{adj}N\right) \left( \limfunc{%
adj}M\right)
\end{equation*}%
and%
\begin{equation*}
M\left( \limfunc{adj}M\right) =\left( \limfunc{adj}M\right) M=\left( \det
M\right) I
\end{equation*}%
are immediate, where $I$ the $3\times 3$ identity matrix, and in fact they
hold more generally for arbitrary $3\times 3$ matrices $M$ and $N$. In the
invertible case we have in addition%
\begin{eqnarray*}
\limfunc{adj}\left( \limfunc{adj}M\right) &=&\det \left( \limfunc{adj}%
M\right) \left( \limfunc{adj}M\right) ^{-1} \\
&=&\det \left( \left( \det M\right) M^{-1}\right) \left( \left( \det
M\right) M^{-1}\right) ^{-1} \\
&=&\left( \det M\right) ^{3}\left( \det M^{-1}\right) \left( \det M\right)
^{-1}M \\
&=&\left( \det M\right) M.
\end{eqnarray*}%
For the fixed symmetric matrix $B$ from (\ref{Symmetric Matrix}), we write%
\begin{equation*}
\limfunc{adj}B=%
\begin{pmatrix}
a_{2}a_{3}-b_{1}^{2} & b_{1}b_{2}-a_{3}b_{3} & b_{1}b_{3}-a_{2}b_{2} \\ 
b_{1}b_{2}-a_{3}b_{3} & a_{1}a_{3}-b_{2}^{2} & b_{2}b_{3}-a_{1}b_{1} \\ 
b_{1}b_{3}-a_{2}b_{2} & b_{2}b_{3}-a_{1}b_{1} & a_{1}a_{2}-b_{3}^{2}%
\end{pmatrix}%
\equiv 
\begin{pmatrix}
\alpha _{1} & \beta _{3} & \beta _{2} \\ 
\beta _{3} & \alpha _{2} & \beta _{1} \\ 
\beta _{2} & \beta _{1} & \alpha _{3}%
\end{pmatrix}%
.
\end{equation*}

Now define the $B$-\textbf{vector product} of vectors $v_{1}$ and $v_{2}$ to
be the vector%
\begin{equation*}
v_{1}\times _{B}v_{2}\equiv \left( v_{1}\times v_{2}\right) \limfunc{adj}B.
\end{equation*}%
The motivation for this definition is given by the following theorem. A
similar result has been obtained in \cite{Collomb}.

\begin{theorem}[Adjugate vector product theorem]
Let $v_{1}$, $v_{2}$ and $v_{3}$ be vectors in $\mathbb{V}^{3}$, and let $%
M\equiv 
\begin{pmatrix}
v_{1} \\ 
v_{2} \\ 
v_{3}%
\end{pmatrix}%
$be the matrix with these vectors as rows. Then for any $3\times 3$
invertible symmetric matrix $B,$ 
\begin{equation*}
\limfunc{adj}\left( MB\right) \equiv 
\begin{pmatrix}
v_{2}\times _{B}v_{3} \\ 
v_{3}\times _{B}v_{1} \\ 
v_{1}\times _{B}v_{2}%
\end{pmatrix}%
^{T}.
\end{equation*}
\end{theorem}

\begin{proof}
By the definition of adjugate matrix, $\func{adj}M$ is%
\begin{equation*}
\limfunc{adj}M=%
\begin{pmatrix}
v_{2}\times v_{3} \\ 
v_{3}\times v_{1} \\ 
v_{1}\times v_{2}%
\end{pmatrix}%
^{T}.
\end{equation*}%
$\allowbreak \allowbreak $Since $\limfunc{adj}\left( MB\right) =\limfunc{adj}%
B\limfunc{adj}M$ and $B$ is symmetric,%
\begin{eqnarray*}
\left( \limfunc{adj}\left( MB\right) \right) ^{T} &=&\left( \limfunc{adj}%
M\right) ^{T}\limfunc{adj}B=%
\begin{pmatrix}
v_{2}\times v_{3} \\ 
v_{3}\times v_{1} \\ 
v_{1}\times v_{2}%
\end{pmatrix}%
\limfunc{adj}B \\
&=&%
\begin{pmatrix}
\left( v_{2}\times v_{3}\right) \limfunc{adj}B \\ 
\left( v_{3}\times v_{1}\right) \limfunc{adj}B \\ 
\left( v_{1}\times v_{2}\right) \limfunc{adj}B%
\end{pmatrix}%
=%
\begin{pmatrix}
v_{2}\times _{B}v_{3} \\ 
v_{3}\times _{B}v_{1} \\ 
v_{1}\times _{B}v_{2}%
\end{pmatrix}%
.
\end{eqnarray*}%
Now take the matrix transpose on both sides.
\end{proof}

The usual linearity and anti-symmetric properties of the Euclidean vector
product hold for $B$-vector products.

\subsection{The $B$-scalar triple product}

The \textbf{Euclidean scalar triple product} of three vectors $v_{1}$, $%
v_{2} $ and $v_{3}$ in $\mathbb{V}^{3}$ (see \cite[pp. 68-71]{Gibbs}) is%
\begin{equation*}
\left[ v_{1},v_{2},v_{3}\right] \equiv v_{1}\cdot \left( v_{2}\times
v_{3}\right) =\det 
\begin{pmatrix}
v_{1} \\ 
v_{2} \\ 
v_{3}%
\end{pmatrix}%
.
\end{equation*}%
We can generalise this definition for an arbitrary symmetric bilinear form
with matrix representation $B$; so we define the $B$\textbf{-scalar triple
product} of $v_{1},v_{2}$ and $v_{3}$ to be 
\begin{equation*}
\left[ v_{1},v_{2},v_{3}\right] _{B}\equiv v_{1}\cdot _{B}\left( v_{2}\times
_{B}v_{3}\right) .
\end{equation*}%
The following result allows for the evaluation of the $B$-scalar triple
product in terms of determinants.

\begin{theorem}[Scalar triple product theorem]
Let $M\equiv 
\begin{pmatrix}
v_{1} \\ 
v_{2} \\ 
v_{3}%
\end{pmatrix}%
$ for vectors $v_{1}$, $v_{2}$ and $v_{3}$ in $\mathbb{V}^{3}$. Then 
\begin{equation*}
\left[ v_{1},v_{2},v_{3}\right] _{B}=\left( \det B\right) \left( \det
M\right) .
\end{equation*}
\end{theorem}

\begin{proof}
From the definitions of the $B$-scalar product, $B$-vector product and the $%
B $-scalar triple product,%
\begin{eqnarray*}
\left[ v_{1},v_{2},v_{3}\right] _{B} &=&v_{1}B\left( \left( v_{2}\times
v_{3}\right) \limfunc{adj}B\right) ^{T} \\
&=&v_{1}\left( B\limfunc{adj}B\right) \left( v_{2}\times v_{3}\right) ^{T}.
\end{eqnarray*}%
As $\limfunc{adj}B=\left( \det B\right) B^{-1}$ and $v_{1}\cdot \left(
v_{2}\times v_{3}\right) =\det M$,%
\begin{eqnarray*}
\left[ v_{1},v_{2},v_{3}\right] _{B} &=&\left( \det B\right) v_{1}\left(
v_{2}\times v_{3}\right) ^{T} \\
&=&\left( \det B\right) \left( v_{1}\cdot \left( v_{2}\times v_{3}\right)
\right) \\
&=&\left( \det B\right) \left( \det M\right)
\end{eqnarray*}%
as required.
\end{proof}

We can now relate $B$-vector products to $B$-perpendicularity.

\begin{corollary}
The vectors $v$ and $w$ in $\mathbb{V}^{3}$ are both $B$-perpendicular to $%
v\times _{B}w$, i.e.%
\begin{equation*}
v\perp _{B}\left( v\times _{B}w\right) \quad \mathrm{and}\text{\quad }w\perp
_{B}\left( v\times _{B}w\right) .
\end{equation*}
\end{corollary}

\begin{proof}
By the Scalar triple product theorem,%
\begin{equation*}
v\cdot _{B}\left( v\times _{B}w\right) =\left[ v,v,w\right] _{B}=\left( \det
B\right) \det 
\begin{pmatrix}
v \\ 
v \\ 
w%
\end{pmatrix}%
=0.
\end{equation*}%
Similarly, $\left[ w,v,w\right] _{B}=0$ and so both $v\perp _{B}\left(
v\times _{B}w\right) $ and $w\perp _{B}\left( v\times _{B}w\right) $.
\end{proof}

We could also rearrange the ordering of $B$-scalar triple products as
follows.

\begin{corollary}
For vectors $v_{1}$, $v_{2}$ and $v_{3}$ in $\mathbb{V}^{3}$, 
\begin{gather*}
\left[ v_{1},v_{2},v_{3}\right] _{B}=\left[ v_{2},v_{3},v_{1}\right] _{B}=%
\left[ v_{3},v_{1},v_{2}\right] _{B} \\
=-\left[ v_{1},v_{3},v_{2}\right] _{B}=-\left[ v_{2},v_{1},v_{3}\right]
_{B}=-\left[ v_{3},v_{2},v_{1}\right] _{B}.
\end{gather*}
\end{corollary}

\begin{proof}
This follows from the corresponding relations for $\left[ v_{1},v_{2},v_{3}%
\right] ,$ or equivalently the transformation properties of the determinant
upon permutation of rows.
\end{proof}

\subsection{The $B$-vector triple product}

Recall that the \textbf{Euclidean vector triple product} of vectors $v_{1}$, 
$v_{2}$ and $v_{3}$ in $\mathbb{V}^{3}$\ (see \cite[pp. 71-75]{Gibbs}) is 
\begin{equation*}
\left\langle v_{1},v_{2},v_{3}\right\rangle \equiv v_{1}\times \left(
v_{2}\times v_{3}\right) .
\end{equation*}%
The $B$-\textbf{vector triple product} of the vectors is similarly defined by%
\begin{equation*}
\left\langle v_{1},v_{2},v_{3}\right\rangle _{B}\equiv v_{1}\times
_{B}\left( v_{2}\times _{B}v_{3}\right) .
\end{equation*}%
We can evaluate this by generalising a classical result of Lagrange \cite%
{Lagrange} from the Euclidean vector triple product to $B$-vector triple
products, following the general lines of argument of \cite{CM} and \cite[pp.
28-29]{Spiegel} in the Euclidean case; the proof is surprisingly complicated.

\begin{theorem}[Lagrange's formula]
For vectors $v_{1}$, $v_{2}$ and $v_{3}$ in $\mathbb{V}^{3}$,%
\begin{equation*}
\left\langle v_{1},v_{2},v_{3}\right\rangle _{B}=\left( \det B\right) \left[
\left( v_{1}\cdot _{B}v_{3}\right) v_{2}-\left( v_{1}\cdot _{B}v_{2}\right)
v_{3}\right] .
\end{equation*}
\end{theorem}

\begin{proof}
Let $w\equiv \left\langle v_{1},v_{2},v_{3}\right\rangle _{B}$. If $v_{2}$
and $v_{3}$ are linearly dependent, then $v_{2}\times _{B}v_{3}=\mathbf{0}$
and thus $\left\langle v_{1},v_{2},v_{3}\right\rangle _{B}=\mathbf{0}$.
Furthermore, we are able to write one of them as a scalar multiple of the
other, which implies that 
\begin{equation*}
\left( v_{1}\cdot _{B}v_{3}\right) v_{2}-\left( v_{1}\cdot _{B}v_{2}\right)
v_{3}=\mathbf{0}
\end{equation*}%
and thus the required result holds. So we may suppose that $v_{2}$ and $%
v_{3} $ are linearly independent. From Corollary 4, $\left( v_{2}\times
_{B}v_{3}\right) \perp _{B}w$ and thus%
\begin{equation*}
v_{2}\perp _{B}\left( v_{2}\times _{B}v_{3}\right) \quad \text{\textrm{and}%
\quad }v_{3}\perp _{B}\left( v_{2}\times _{B}v_{3}\right) .
\end{equation*}%
As $w$ is parallel to $v_{2}$ and $v_{3}$, we can deduce that $w$ is equal
to some linear combination of $v_{2}$ and $v_{3}$. So, for some scalars $%
\alpha $ and $\beta $ in $\mathbb{F}$, we have 
\begin{equation*}
w=\alpha v_{2}+\beta v_{3}.
\end{equation*}%
Furthermore, since $v_{1}\perp _{B}w$, the definition of $B$%
-perpendicularity implies that 
\begin{equation*}
w\cdot _{B}v_{1}=\alpha \left( v_{1}\cdot _{B}v_{2}\right) +\beta \left(
v_{1}\cdot _{B}v_{3}\right) =0.
\end{equation*}%
This equality is true precisely when $\alpha =\lambda \left( v_{1}\cdot
_{B}v_{3}\right) $ and $\beta =-\lambda \left( v_{1}\cdot _{B}v_{2}\right) $%
, for some non-zero scalar $\lambda $ in $\mathbb{F}$. Hence,%
\begin{equation*}
w=\lambda \left[ \left( v_{1}\cdot _{B}v_{3}\right) v_{2}-\left( v_{1}\cdot
_{B}v_{2}\right) v_{3}\right] .
\end{equation*}%
To proceed, we first want to prove that $\lambda $ is independent of the
choices $v_{1}$, $v_{2}$ and $v_{3}$, so that we can compute $w$ for
arbitrary $v_{1}$, $v_{2}$ and $v_{3}$. First, suppose that $\lambda $ is
dependent on $v_{1}$, $v_{2}$ and $v_{3}$, so that we may define $\lambda
\equiv \lambda \left( v_{1},v_{2},v_{3}\right) $. Given another vector $d$
in $\mathbb{V}^{3}$, we have%
\begin{equation*}
w\cdot _{B}d=\lambda \left( v_{1},v_{2},v_{3}\right) \left[ \left(
v_{1}\cdot _{B}v_{3}\right) \left( v_{2}\cdot _{B}d\right) -\left(
v_{1}\cdot _{B}v_{2}\right) \left( v_{3}\cdot _{B}d\right) \right] .
\end{equation*}%
Directly substituting the definition of $w$, we use the Scalar triple
product theorem to obtain%
\begin{equation*}
w\cdot _{B}d=\left( v_{1}\times _{B}\left( v_{2}\times _{B}v_{3}\right)
\right) \cdot _{B}d=v_{1}\cdot _{B}\left( \left( v_{2}\times
_{B}v_{3}\right) \times _{B}d\right) =-v_{1}\cdot _{B}\left\langle
d,v_{2},v_{3}\right\rangle _{B}.
\end{equation*}%
Based on our calculations of $w$, we then deduce that%
\begin{eqnarray*}
-v_{1}\cdot _{B}\left\langle d,v_{2},v_{3}\right\rangle _{B} &=&-v_{1}\cdot
_{B}\left( \lambda \left( d,v_{2},v_{3}\right) \left[ \left( d\cdot
_{B}v_{3}\right) v_{2}-\left( d\cdot _{B}v_{2}\right) v_{3}\right] \right) \\
&=&\lambda \left( d,v_{2},v_{3}\right) \left[ \left( v_{1}\cdot
_{B}v_{3}\right) \left( v_{2}\cdot _{B}d\right) -\left( v_{1}\cdot
_{B}v_{2}\right) \left( v_{3}\cdot _{B}d\right) \right] .
\end{eqnarray*}%
Because this expression is equal to $w\cdot _{B}d$, we deduce that $\lambda
\left( v_{1},v_{2},v_{3}\right) =\lambda \left( d,v_{2},v_{3}\right) $ and
hence $\lambda $ must be independent of the choice of $v_{1}$. With this,
now suppose instead that $\lambda \equiv \lambda \left( v_{2},v_{3}\right) $%
, so that%
\begin{equation*}
w\cdot _{B}d=\lambda \left( v_{2},v_{3}\right) \left[ \left( v_{1}\cdot
_{B}v_{3}\right) \left( v_{2}\cdot _{B}d\right) -\left( v_{1}\cdot
_{B}v_{2}\right) \left( v_{3}\cdot _{B}d\right) \right]
\end{equation*}%
for a vector $d$ in $\mathbb{V}^{3}$. By direct substitution of $w$, we use
the Scalar triple product theorem to obtain%
\begin{equation*}
w\cdot _{B}d=\left( v_{1}\times _{B}\left( v_{2}\times _{B}v_{3}\right)
\right) \cdot _{B}d=\left( v_{2}\times _{B}v_{3}\right) \cdot _{B}\left(
d\times _{B}v_{1}\right) =v_{2}\cdot _{B}\left\langle
v_{3},d,v_{1}\right\rangle _{B}.
\end{equation*}%
Similarly, based on the calculations of $w$ previously, we have%
\begin{eqnarray*}
v_{2}\cdot _{B}\left\langle v_{3},d,v_{1}\right\rangle _{B} &=&v_{2}\cdot
_{B}\lambda \left( d_{2},v_{3}\right) \left( \left( v_{1}\cdot
_{B}v_{3}\right) d-\left( v_{3}\cdot _{B}d\right) v_{1}\right) \\
&=&\lambda \left( d,v_{1}\right) \left[ \left( v_{1}\cdot _{B}v_{3}\right)
\left( v_{2}\cdot _{B}d\right) -\left( v_{1}\cdot _{B}v_{2}\right) \left(
v_{3}\cdot _{B}d\right) \right] .
\end{eqnarray*}%
Because this expression is also equal to $w\cdot _{B}d$, we deduce that $%
\lambda \left( v_{2},v_{3}\right) =\lambda \left( d,v_{1}\right) $ and
conclude that $\lambda $ is indeed independent of $v_{2}$ and $v_{3}$, in
addition to $v_{1}$. So, we substitute any choice of vectors $v_{1}$, $v_{2}$
and $v_{3}$ in order to find $\lambda $. With this, suppose that $%
v_{2}\equiv \left( 1,0,0\right) $ and $v_{1}=v_{3}\equiv \left( 0,1,0\right) 
$. Then,%
\begin{eqnarray*}
v_{2}\times _{B}v_{3} &=&\left[ \left( 1,0,0\right) \times \left(
0,1,0\right) \right] \func{adj}B \\
&=&\left( 0,0,1\right) 
\begin{pmatrix}
\alpha _{1} & \beta _{3} & \beta _{2} \\ 
\beta _{3} & \alpha _{2} & \beta _{1} \\ 
\beta _{2} & \beta _{1} & \alpha _{3}%
\end{pmatrix}
\\
&=&\left( \beta _{2},\beta _{1},\alpha _{3}\right)
\end{eqnarray*}%
and hence%
\begin{eqnarray*}
\left\langle v_{1},v_{2},v_{3}\right\rangle _{B} &=&\left[ \left(
0,1,0\right) \times \left( \beta _{2},\beta _{1},\alpha _{3}\right) \right] 
\begin{pmatrix}
\alpha _{1} & \beta _{3} & \beta _{2} \\ 
\beta _{3} & \alpha _{2} & \beta _{1} \\ 
\beta _{2} & \beta _{1} & \alpha _{3}%
\end{pmatrix}
\\
&=&\left( \alpha _{3},0,-\beta _{2}\right) 
\begin{pmatrix}
\alpha _{1} & \beta _{3} & \beta _{2} \\ 
\beta _{3} & \alpha _{2} & \beta _{1} \\ 
\beta _{2} & \beta _{1} & \alpha _{3}%
\end{pmatrix}
\\
&=&\left( \alpha _{1}\alpha _{3}-\beta _{2}^{2},\alpha _{3}\beta _{3}-\beta
_{1}\beta _{2},0\right) .
\end{eqnarray*}%
Now use the fact that $\limfunc{adj}\left( \limfunc{adj}B\right) =\left(
\det B\right) B$ to obtain%
\begin{eqnarray*}
\func{adj}%
\begin{pmatrix}
\alpha _{1} & \beta _{3} & \beta _{2} \\ 
\beta _{3} & \alpha _{2} & \beta _{1} \\ 
\beta _{2} & \beta _{1} & \alpha _{3}%
\end{pmatrix}
&=&%
\begin{pmatrix}
\alpha _{2}\alpha _{3}-\beta _{1}^{2} & \beta _{1}\beta _{2}-\alpha
_{3}\beta _{3} & \beta _{1}\beta _{3}-\alpha _{2}\beta _{2} \\ 
\beta _{1}\beta _{2}-\alpha _{3}\beta _{3} & \alpha _{1}\alpha _{3}-\beta
_{2}^{2} & \beta _{2}\beta _{3}-\alpha _{1}\beta _{1} \\ 
\beta _{1}\beta _{3}-\alpha _{2}\beta _{2} & \beta _{2}\beta _{3}-\alpha
_{1}\beta _{1} & \alpha _{1}\alpha _{2}-\beta _{3}^{2}%
\end{pmatrix}%
\allowbreak \\
&=&\left( \det B\right) 
\begin{pmatrix}
a_{1} & b_{3} & b_{2} \\ 
b_{3} & a_{2} & b_{1} \\ 
b_{2} & b_{1} & a_{3}%
\end{pmatrix}%
\end{eqnarray*}%
so that 
\begin{equation*}
\left\langle v_{1},v_{2},v_{3}\right\rangle _{B}=\left( \det B\right) \left(
a_{2},-b_{3},0\right) .
\end{equation*}%
Since $v_{1}\cdot _{B}v_{2}=e_{1}Be_{2}^{T}=b_{3}$ and $v_{1}\cdot
_{B}v_{3}=e_{2}Be_{2}^{T}=a_{2}$, it follows that%
\begin{eqnarray*}
\left( \det B\right) \left( a_{2},-b_{3},0\right) &=&\left( \det B\right) 
\left[ \left( v_{1}\cdot _{B}v_{3}\right) e_{1}-\left( v_{1}\cdot
_{B}v_{2}\right) e_{2}\right] \\
&=&\left( \det B\right) \left[ \left( v_{1}\cdot _{B}v_{3}\right)
v_{2}-\left( v_{1}\cdot _{B}v_{2}\right) v_{3}\right] .
\end{eqnarray*}%
From this, we deduce that $\lambda =\det B$ and hence%
\begin{equation*}
\left\langle v_{1},v_{2},v_{3}\right\rangle _{B}=\left( \det B\right) \left[
\left( v_{1}\cdot _{B}v_{3}\right) v_{2}-\left( v_{1}\cdot _{B}v_{2}\right)
v_{3}\right]
\end{equation*}%
as required.
\end{proof}

The $B$-vector product is not an associative operation, but by the
anti-symmetric property of $B$-vector products, we see that%
\begin{equation*}
\left\langle v_{1},v_{2},v_{3}\right\rangle _{B}=-\left\langle
v_{1},v_{3},v_{2}\right\rangle _{B}.
\end{equation*}%
The following result, attributed in the Euclidean case to Jacobi \cite%
{Jacobi}, connects the theory of $B$-vector products to the theory of Lie
algebras and links the three $B$-vector triple products which differ by an
even permutation of the indices.

\begin{theorem}[Jacobi identity]
For vectors $v_{1}$, $v_{2}$ and $v_{3}$ in $\mathbb{V}^{3}$,%
\begin{equation*}
\left\langle v_{1},v_{2},v_{3}\right\rangle _{B}+\left\langle
v_{2},v_{3},v_{1}\right\rangle _{B}+\left\langle
v_{3},v_{1},v_{2}\right\rangle _{B}=\mathbf{0}.
\end{equation*}
\end{theorem}

\begin{proof}
Apply Lagrange's formula to each of the three summands to get%
\begin{equation*}
\left\langle v_{1},v_{2},v_{3}\right\rangle _{B}=\left( \det B\right) \left[
\left( v_{1}\cdot _{B}v_{3}\right) v_{2}-\left( v_{1}\cdot _{B}v_{2}\right)
v_{3}\right]
\end{equation*}%
as well as%
\begin{equation*}
\left\langle v_{2},v_{3},v_{1}\right\rangle _{B}=\left( \det B\right) \left[
\left( v_{1}\cdot _{B}v_{2}\right) v_{3}-\left( v_{2}\cdot _{B}v_{3}\right)
v_{1}\right]
\end{equation*}%
and%
\begin{equation*}
\left\langle v_{3},v_{1},v_{2}\right\rangle _{B}=\left( \det B\right) \left[
\left( v_{2}\cdot _{B}v_{3}\right) v_{1}-\left( v_{1}\cdot _{B}v_{3}\right)
v_{2}\right] .
\end{equation*}%
So,%
\begin{eqnarray*}
&&\left\langle v_{1},v_{2},v_{3}\right\rangle _{B}+\left\langle
v_{2},v_{3},v_{1}\right\rangle _{B}+\left\langle
v_{3},v_{1},v_{2}\right\rangle _{B} \\
&=&\left( \det B\right) \left[ \left( v_{1}\cdot _{B}v_{3}\right)
v_{2}-\left( v_{1}\cdot _{B}v_{2}\right) v_{3}\right] +\left( \det B\right) %
\left[ \left( v_{1}\cdot _{B}v_{2}\right) v_{3}-\left( v_{2}\cdot
_{B}v_{3}\right) v_{1}\right] \\
&&+\left( \det B\right) \left[ \left( v_{2}\cdot _{B}v_{3}\right)
v_{1}-\left( v_{1}\cdot _{B}v_{3}\right) v_{2}\right] =\mathbf{0}
\end{eqnarray*}%
as required.
\end{proof}

\subsection{The $B$-scalar quadruple product}

Recall that the \textbf{Euclidean scalar quadruple product} of vectors $%
v_{1} $, $v_{2}$,$\,v_{3}$ and $v_{4}$ in $\mathbb{V}^{3}$ (see \cite[pp.
75-76]{Gibbs}) is the scalar 
\begin{equation*}
\left[ v_{1},v_{2};v_{3},v_{4}\right] \equiv \left( v_{1}\times v_{2}\right)
\cdot \left( v_{3}\times v_{4}\right) .
\end{equation*}%
We define similarly the $B$\textbf{-scalar quadruple product} to be the
quantity 
\begin{equation*}
\left[ v_{1},v_{2};v_{3},v_{4}\right] _{B}\equiv \left( v_{1}\times
_{B}v_{2}\right) \cdot _{B}\left( v_{3}\times _{B}v_{4}\right) .
\end{equation*}%
The following result, which originated from separate works of Binet \cite%
{Binet} and Cauchy \cite{Cauchy} in the Euclidean setting (also see \cite{BS}
and \cite[p. 29]{Spiegel}), allows us to compute $B$-scalar quadruple
products purely in terms of $B$-scalar products.

\begin{theorem}[Binet-Cauchy identity]
For vectors $v_{1}$, $v_{2}$, $v_{3}$ and $v_{4}$ in $\mathbb{V}^{3}$,%
\begin{equation*}
\left[ v_{1},v_{2};v_{3},v_{4}\right] _{B}=\left( \det B\right) \left[
\left( v_{1}\cdot _{B}v_{3}\right) \left( v_{2}\cdot _{B}v_{4}\right)
-\left( v_{1}\cdot _{B}v_{4}\right) \left( v_{2}\cdot _{B}v_{3}\right) %
\right] .
\end{equation*}
\end{theorem}

\begin{proof}
Let $w\equiv v_{1}\times _{B}v_{2}$, so that by the Scalar triple product
theorem and Corollary 3,%
\begin{equation*}
\left[ v_{1},v_{2};v_{3},v_{4}\right] _{B}=\left[ w,v_{3},v_{4}\right] _{B}=%
\left[ v_{4},w,v_{3}\right] _{B}.
\end{equation*}%
By Lagrange's formula,%
\begin{equation*}
w\times _{B}v_{3}=-\left\langle v_{3},v_{1},v_{2}\right\rangle _{B}=\left(
\det B\right) \left[ \left( v_{1}\cdot _{B}v_{3}\right) v_{2}-\left(
v_{2}\cdot _{B}v_{3}\right) v_{1}\right]
\end{equation*}%
and hence%
\begin{eqnarray*}
\left[ v_{1},v_{2};v_{3},v_{4}\right] _{B} &=&\left( \left( \det B\right) %
\left[ \left( v_{1}\cdot _{B}v_{3}\right) v_{2}-\left( v_{2}\cdot
_{B}v_{3}\right) v_{1}\right] \right) \cdot v_{4} \\
&=&\left( \det B\right) \left[ \left( v_{1}\cdot _{B}v_{3}\right) \left(
v_{2}\cdot _{B}v_{4}\right) -\left( v_{1}\cdot _{B}v_{4}\right) \left(
v_{2}\cdot _{B}v_{3}\right) \right]
\end{eqnarray*}%
as required.
\end{proof}

An important special case of the classical Binet-Cauchy identity is a result
of Lagrange \cite{Lagrange}, which we now generalise. We distinguish this
from Lagrange's formula, which computes the $B$-vector triple product of
three vectors, by calling it Lagrange's identity.

\begin{theorem}[Lagrange's identity]
Given vectors $v_{1}$ and $v_{2}$ in $\mathbb{V}^{3}$,%
\begin{equation*}
Q_{B}\left( v_{1}\times _{B}v_{2}\right) =\left( \det B\right) \left[
Q_{B}\left( v_{1}\right) Q_{B}\left( v_{2}\right) -\left( v_{1}\cdot
_{B}v_{2}\right) ^{2}\right] .
\end{equation*}
\end{theorem}

\begin{proof}
This immediately follows from the Binet-Cauchy identity by setting $%
v_{1}=v_{3}$ and $v_{2}=v_{4}$.
\end{proof}

Here is another consequence of the Binet-Cauchy identity, which is somewhat
similar to the Jacobi identity for $B$-vector triple products.

\begin{corollary}
For vectors $v_{1}$, $v_{2}$, $v_{3}$ and $v_{4}$ in $\mathbb{V}^{3}$, we
have%
\begin{equation*}
\left[ v_{1},v_{2};v_{3},v_{4}\right] _{B}+\left[ v_{2},v_{3};v_{1},v_{4}%
\right] _{B}+\left[ v_{3},v_{1};v_{2},v_{4}\right] _{B}=0.
\end{equation*}
\end{corollary}

\begin{proof}
From the Binet-Cauchy identity, the summands evaluate to 
\begin{equation*}
\left[ v_{1},v_{2};v_{3},v_{4}\right] _{B}=\left( \det B\right) \left[
\left( v_{1}\cdot _{B}v_{3}\right) \left( v_{2}\cdot _{B}v_{4}\right)
-\left( v_{1}\cdot _{B}v_{4}\right) \left( v_{2}\cdot _{B}v_{3}\right) %
\right]
\end{equation*}%
as well as%
\begin{equation*}
\left[ v_{2},v_{3};v_{1},v_{4}\right] _{B}=\left( \det B\right) \left[
\left( v_{2}\cdot _{B}v_{1}\right) \left( v_{3}\cdot _{B}v_{4}\right)
-\left( v_{2}\cdot _{B}v_{4}\right) \left( v_{3}\cdot _{B}v_{1}\right) %
\right]
\end{equation*}%
and%
\begin{equation*}
\left[ v_{3},v_{1};v_{2},v_{4}\right] _{B}=\left( \det B\right) \left[
\left( v_{3}\cdot _{B}v_{2}\right) \left( v_{1}\cdot _{B}v_{4}\right)
-\left( v_{3}\cdot _{B}v_{4}\right) \left( v_{1}\cdot _{B}v_{2}\right) %
\right] .
\end{equation*}%
When we add these three quantities, we get $0$ as required.
\end{proof}

\subsection{The $B$-vector quadruple product}

Recall that the \textbf{Euclidean vector quadruple product} of vectors $%
v_{1} $, $v_{2}$,$\,v_{3}$ and $v_{4}$ in $\mathbb{V}^{3}$ (see \cite[pp.
76-77]{Gibbs}) is the vector 
\begin{equation*}
\left\langle v_{1},v_{2};v_{3},v_{4}\right\rangle \equiv \left( v_{1}\times
v_{2}\right) \times \left( v_{3}\times v_{4}\right) .
\end{equation*}%
Define similarly the $B$-\textbf{vector quadruple product} to be 
\begin{equation*}
\left\langle v_{1},v_{2};v_{3},v_{4}\right\rangle _{B}\equiv \left(
v_{1}\times _{B}v_{2}\right) \times _{B}\left( v_{3}\times _{B}v_{4}\right) .
\end{equation*}%
The key result here is the following.

\begin{theorem}[$B$-vector quadruple product theorem]
For vectors $v_{1}$, $v_{2}$, $v_{3}$ and $v_{4}$ in $\mathbb{V}^{3}$,%
\begin{eqnarray*}
\left\langle v_{1},v_{2};v_{3},v_{4}\right\rangle _{B} &=&\left( \det
B\right) \left( \left[ v_{1},v_{2},v_{4}\right] _{B}v_{3}-\left[
v_{1},v_{2},v_{3}\right] _{B}v_{4}\right) \\
&=&\left( \det B\right) \left( \left[ v_{1},v_{3},v_{4}\right] _{B}v_{2}-%
\left[ v_{2},v_{3},v_{4}\right] _{B}v_{1}\right) .
\end{eqnarray*}
\end{theorem}

\begin{proof}
If $u\equiv v_{1}\times _{B}v_{2}$, then use Lagrange's formula to get%
\begin{equation*}
\left\langle v_{1},v_{2};v_{3},v_{4}\right\rangle _{B}=\left\langle
u,v_{3},v_{4}\right\rangle _{B}=\left( \det B\right) \left[ \left( u\cdot
_{B}v_{4}\right) v_{3}-\left( u\cdot _{B}v_{3}\right) v_{4}\right] .
\end{equation*}%
From the Scalar triple product theorem, 
\begin{equation*}
u\cdot _{B}v_{3}=\left[ v_{1},v_{2},v_{3}\right] _{B}\quad \mathrm{and}\text{%
\quad }u\cdot _{B}v_{4}=\left[ v_{1},v_{2},v_{4}\right] _{B}.
\end{equation*}%
Therefore,%
\begin{equation*}
\left\langle v_{1},v_{2};v_{3},v_{4}\right\rangle _{B}=\left( \det B\right) 
\left[ \left[ v_{1},v_{2},v_{4}\right] _{B}v_{3}-\left[ v_{1},v_{2},v_{3}%
\right] _{B}v_{4}\right] .
\end{equation*}%
Since 
\begin{equation*}
\left\langle v_{1},v_{2};v_{3},v_{4}\right\rangle _{B}=-\left\langle
v_{3},v_{4};v_{1},v_{2}\right\rangle _{B}=\left( v_{3}\times
_{B}v_{4}\right) \times _{B}\left( v_{2}\times _{B}v_{1}\right)
\end{equation*}%
Corollary 3 gives us%
\begin{eqnarray*}
\left\langle v_{1},v_{2};v_{3},v_{4}\right\rangle _{B} &=&\left( \det
B\right) \left( \left[ v_{3},v_{4},v_{1}\right] _{B}v_{2}-\left[
v_{3},v_{4},v_{2}\right] _{B}v_{1}\right) \\
&=&\left( \det B\right) \left( \left[ v_{1},v_{3},v_{4}\right] _{B}v_{2}-%
\left[ v_{2},v_{3},v_{4}\right] _{B}v_{1}\right) .
\end{eqnarray*}
\end{proof}

As a corollary, we find a relation satisfied by any four vectors in
three-dimensional vector space, extending the result in \cite[p. 76]{Gibbs}
to a general metrical framework.

\begin{corollary}[Four vector relation]
Suppose that $v_{1}$, $v_{2}$, $v_{3}$ and $v_{4}$ are vectors in $\mathbb{V}%
^{3}.$ Then 
\begin{equation*}
\left[ v_{2},v_{3},v_{4}\right] _{B}v_{1}-\left[ v_{1},v_{3},v_{4}\right]
_{B}v_{2}+\left[ v_{1},v_{2},v_{4}\right] _{B}v_{3}-\left[ v_{1},v_{2},v_{3}%
\right] _{B}v_{4}=\mathbf{0.}
\end{equation*}
\end{corollary}

\begin{proof}
This is an immediate consequence of equating the two equations from the
previous result, after cancelling the non-zero factor $\det B$.
\end{proof}

As another consequence, we get an expression for the meet of two distinct
two-dimensional subspaces.

\begin{corollary}
If $U\equiv \limfunc{span}\left( v_{1},v_{2}\right) $ and $V\equiv \limfunc{%
span}\left( v_{3},v_{4}\right) $ are distinct two-dimensional subspaces then 
$v\equiv \left\langle v_{1},v_{2};v_{3},v_{4}\right\rangle _{B}$ spans $%
U\cap V$.
\end{corollary}

\begin{proof}
Clearly $v$ is both in $U$ and in $V$ from the $B$-vector quadruple product
theorem. We need only show that it is non-zero, but this follows from 
\begin{equation*}
\left\langle v_{1},v_{2};v_{3},v_{4}\right\rangle _{B}=\left( \det B\right)
\left( \left[ v_{1},v_{2},v_{4}\right] _{B}v_{3}-\left[ v_{1},v_{2},v_{3}%
\right] _{B}v_{4}\right)
\end{equation*}%
since by assumption $v_{3}$ and $v_{4}$ are linearly independent, and at
least one of $\left[ v_{1},v_{2},v_{4}\right] _{B}$ and $\left[
v_{1},v_{2},v_{3}\right] _{B}$ must be non-zero since otherwise both $v_{4}$
and $v_{3}$ lie in $U$, which contradicts the assumption that the $U$ and $V$
are distinct.
\end{proof}

A special case occurs when each of the factors of the $B$-quadruple vector
product contains a common vector. This extends the result in \cite[p. 80]%
{Gibbs} to a general metrical framework.

\begin{corollary}
If $v_{1}$, $v_{2}$ and $v_{3}$ are vectors in $\mathbb{V}^{3}$ then%
\begin{equation*}
\left\langle v_{1},v_{2};v_{1},v_{3}\right\rangle _{B}=\left( \det B\right) 
\left[ v_{1},v_{2},v_{3}\right] _{B}v_{1}.
\end{equation*}
\end{corollary}

\begin{proof}
This follows from 
\begin{equation*}
\left\langle v_{1},v_{2};v_{1},v_{3}\right\rangle _{B}=\left( \det B\right)
\left( \left[ v_{1},v_{2},v_{3}\right] _{B}v_{1}-\left[ v_{1},v_{2},v_{1}%
\right] _{B}v_{3}\right)
\end{equation*}%
together with the fact that $\left[ v_{1},v_{2},v_{1}\right] _{B}=0$.
\end{proof}

Yet another consequence is given below, which was alluded to \cite[p. 86]%
{Gibbs} in for the Euclidean case.

\begin{theorem}[Triple scalar product of products]
For vectors $v_{1}$,$v_{2}$ and $v_{3}$ in $\mathbb{V}^{3}$, 
\begin{equation*}
\left[ v_{2}\times _{B}v_{3},v_{3}\times _{B}v_{1},v_{1}\times _{B}v_{2}%
\right] _{B}=\left( \det B\right) \left( \left[ v_{1},v_{2},v_{3}\right]
_{B}\right) ^{2}.
\end{equation*}
\end{theorem}

\begin{proof}
From the previous corollary, 
\begin{eqnarray*}
\left[ v_{2}\times _{B}v_{3},v_{3}\times _{B}v_{1},v_{1}\times _{B}v_{2}%
\right] _{B} &=&-\left( v_{2}\times _{B}v_{3}\right) \cdot _{B}\left( \left(
v_{1}\times _{B}v_{3}\right) \times _{B}\left( v_{1}\times _{B}v_{2}\right)
\right) \\
&=&-\left( v_{2}\times _{B}v_{3}\right) \cdot _{B}\left( \det B\right) \left[
v_{1},v_{3},v_{2}\right] _{B}v_{1} \\
&=&\left( \det B\right) \left( \left[ v_{1},v_{2},v_{3}\right] _{B}\right)
^{2}.
\end{eqnarray*}
\end{proof}

It follows that if $v_{1}$, $v_{2}$ and $v_{3}$ are linearly independent,
then so are $v_{1}\times _{B}v_{2}$, $v_{2}\times _{B}v_{3}$ and $%
v_{3}\times _{B}v_{1}.$ This also suggests there is a kind of duality here,
which we can clarify by the following result, which is a generalization of
Exercise 8 of \cite[p. 116]{Narayan}.

\begin{theorem}
Suppose that $v_{1},v_{2}$ and $v_{3}$ are linearly independent vectors in $%
\mathbb{V}^{3},$ so that $\left[ v_{1},v_{2},v_{3}\right] _{B}$ is non-zero.
Define 
\begin{equation*}
w_{1}=\frac{v_{2}\times _{B}v_{3}}{\left[ v_{1},v_{2},v_{3}\right] _{B}}%
,\quad w_{2}=\frac{v_{3}\times _{B}v_{1}}{\left[ v_{1},v_{2},v_{3}\right]
_{B}}\quad \mathrm{and}\quad w_{3}=\frac{v_{1}\times _{B}v_{2}}{\left[
v_{1},v_{2},v_{3}\right] _{B}}.
\end{equation*}%
Then%
\begin{equation*}
v_{1}\times _{B}w_{1}+v_{2}\times _{B}w_{2}+v_{3}\times _{B}w_{3}=\mathbf{0}
\end{equation*}%
and%
\begin{equation*}
v_{1}\cdot _{B}w_{1}+v_{2}\cdot _{B}w_{2}+v_{3}\cdot _{B}w_{3}=3
\end{equation*}%
and%
\begin{equation*}
\left[ v_{1},v_{2},v_{3}\right] _{B}\left[ w_{1},w_{2},w_{3}\right]
_{B}=\det B
\end{equation*}%
and 
\begin{equation*}
v_{1}=\frac{w_{2}\times w_{3}}{\left[ w_{1},w_{2},w_{3}\right] _{B}},\quad
v_{2}=\frac{w_{3}\times w_{1}}{\left[ w_{1},w_{2},w_{3}\right] _{B}}\quad 
\mathrm{and}\quad v_{3}=\frac{w_{1}\times w_{2}}{\left[ w_{1},w_{2},w_{3}%
\right] _{B}}.
\end{equation*}
\end{theorem}

\begin{proof}
By the Jacobi identity, 
\begin{eqnarray*}
&&v_{1}\times _{B}w_{1}+v_{2}\times _{B}w_{2}+v_{3}\times _{B}w_{3} \\
&=&\frac{1}{\left[ v_{1},v_{2},v_{3}\right] _{B}}\left( v_{1}\times
_{B}\left( v_{2}\times _{B}v_{3}\right) +v_{2}\times _{B}\left( v_{3}\times
_{B}v_{1}\right) +v_{3}\times _{B}\left( v_{1}\times _{B}v_{2}\right) \right)
\\
&=&\frac{1}{\left[ v_{1},v_{2},v_{3}\right] _{B}}\left( \left\langle
v_{1},v_{2},v_{3}\right\rangle _{B}+\left\langle
v_{2},v_{3},v_{1}\right\rangle _{B}+\left\langle
v_{3},v_{1},v_{2}\right\rangle _{B}\right) \\
&=&\mathbf{0.}
\end{eqnarray*}%
Moreover, we use the Scalar triple product to obtain 
\begin{eqnarray*}
&&v_{1}\cdot _{B}w_{1}+v_{2}\cdot _{B}w_{2}+v_{3}\cdot _{B}w_{3} \\
&=&\frac{1}{\left[ v_{1},v_{2},v_{3}\right] _{B}}\left( v_{1}\cdot
_{B}\left( v_{2}\times _{B}v_{3}\right) +v_{2}\cdot _{B}\left( v_{3}\times
_{B}v_{1}\right) +v_{3}\cdot _{B}\left( v_{1}\times _{B}v_{2}\right) \right)
\\
&=&\frac{1}{\left[ v_{1},v_{2},v_{3}\right] _{B}}\left( \left[
v_{1},v_{2},v_{3}\right] _{B}+\left[ v_{2},v_{3},v_{1}\right] _{B}+\left[
v_{3},v_{1},v_{2}\right] _{B}\right) \\
&=&\frac{1}{\left[ v_{1},v_{2},v_{3}\right] _{B}}\left( \left[
v_{1},v_{2},v_{3}\right] _{B}+\left[ v_{2},v_{3},v_{1}\right] _{B}+\left[
v_{3},v_{1},v_{2}\right] _{B}\right) \\
&=&3.
\end{eqnarray*}%
By the Triple scalar product of products theorem, 
\begin{eqnarray*}
\left[ w_{1},w_{2},w_{3}\right] _{B} &=&\left[ \frac{v_{2}\times _{B}v_{3}}{%
\left[ v_{1},v_{2},v_{3}\right] _{B}},\frac{v_{3}\times _{B}v_{1}}{\left[
v_{1},v_{2},v_{3}\right] _{B}},\frac{v_{1}\times _{B}v_{2}}{\left[
v_{1},v_{2},v_{3}\right] _{B}}\right] _{B} \\
&=&\frac{1}{\left[ v_{1},v_{2},v_{3}\right] _{B}^{3}}\left[ v_{2}\times
_{B}v_{3},v_{3}\times _{B}v_{1},v_{1}\times _{B}v_{2}\right] _{B} \\
&=&\frac{1}{\left[ v_{1},v_{2},v_{3}\right] _{B}^{3}}\left( \det B\right)
\left( \left[ v_{1},v_{2},v_{3}\right] _{B}\right) ^{2} \\
&=&\frac{\det B}{\left[ v_{1},v_{2},v_{3}\right] _{B}}.
\end{eqnarray*}%
Hence, 
\begin{equation*}
\left[ v_{1},v_{2},v_{3}\right] _{B}\left[ w_{1},w_{2},w_{3}\right]
_{B}=\det B.
\end{equation*}%
Given this result, we have 
\begin{eqnarray*}
w_{2}\times _{B}w_{3} &=&\frac{v_{3}\times _{B}v_{1}}{\left[
v_{1},v_{2},v_{3}\right] _{B}}\times _{B}\frac{v_{1}\times _{B}v_{2}}{\left[
v_{1},v_{2},v_{3}\right] _{B}} \\
&=&\frac{1}{\left[ v_{1},v_{2},v_{3}\right] _{B}^{2}}\left[ \left(
v_{3}\times _{B}v_{1}\right) \times _{B}\left( v_{1}\times _{B}v_{2}\right) %
\right] \\
&=&-\frac{1}{\left[ v_{1},v_{2},v_{3}\right] _{B}^{2}}\left[ \left(
v_{1}\times _{B}v_{3}\right) \times _{B}\left( v_{1}\times _{B}v_{2}\right) %
\right] \\
&=&-\frac{1}{\left[ v_{1},v_{2},v_{3}\right] _{B}^{2}}\left( \det B\right) %
\left[ v_{1},v_{3},v_{2}\right] _{B}v_{1} \\
&=&\frac{\left( \det B\right) }{\left[ v_{1},v_{2},v_{3}\right] _{B}}v_{1} \\
&=&\left[ w_{1},w_{2},w_{3}\right] _{B}v_{1}
\end{eqnarray*}%
where the fourth line results from Corollary 13 and the fifth line results
from the Scalar triple product theorem. Thus 
\begin{equation*}
v_{1}=\frac{w_{2}\times w_{3}}{\left[ w_{1},w_{2},w_{3}\right] _{B}}
\end{equation*}%
and the results for $v_{2}$ and $v_{3}$ are similar.
\end{proof}

The first and last of the results of \cite[p. 116]{Narayan} was also proven
in \cite[p. 86]{Gibbs} for the Euclidean situation.

\section{Rational trigonometry for a vector triangle}

We would like to extend the framework of \cite{WildDP} and configure
three-dimensional rational trigonometry so it works over a general field and
general bilinear form, centrally framing our discussion on the tools we have
developed above.

We assume as before a $B$-scalar product on the three-dimensional vector
space $\mathbb{V}^{3}$. A \textbf{vector triangle} $\overline{v_{1}v_{2}v_{3}%
}$ is an unordered collection of three vectors $v_{1}$, $v_{2}$ and $v_{3}$
satisfying 
\begin{equation*}
v_{1}+v_{2}+v_{3}=\mathbf{0}.
\end{equation*}%
The $B$-\textbf{quadrances }of such a triangle are the numbers 
\begin{equation*}
Q_{1}\equiv Q_{B}\left( v_{1}\right) ,\quad Q_{2}\equiv Q_{B}\left(
v_{2}\right) \quad \mathrm{and}\quad Q_{3}\equiv Q_{B}\left( v_{3}\right) .
\end{equation*}%
Define \textbf{Archimedes' function} \cite[p. 64]{WildDP} as%
\begin{equation*}
A\left( a,b,c\right) \equiv \left( a+b+c\right) ^{2}-2\left(
a^{2}+b^{2}+c^{2}\right) .
\end{equation*}%
The $B$\textbf{-quadrea} of $\overline{v_{1}v_{2}v_{3}}$ is then%
\begin{equation*}
\mathcal{A}_{B}\left( \overline{v_{1}v_{2}v_{3}}\right) =A\left(
Q_{1},Q_{2},Q_{3}\right) .
\end{equation*}

\begin{theorem}[Quadrea theorem]
Given a vector triangle $\overline{v_{1}v_{2}v_{3}}$ 
\begin{equation*}
Q_{B}\left( v_{1}\times _{B}v_{2}\right) =Q_{B}\left( v_{2}\times
_{B}v_{3}\right) =Q_{B}\left( v_{3}\times _{B}v_{1}\right) =\frac{\left(
\det B\right) }{4}\mathcal{A}_{B}\left( \overline{v_{1}v_{2}v_{3}}\right) .
\end{equation*}
\end{theorem}

\begin{proof}
By bilinearity, we have 
\begin{equation*}
Q_{B}\left( v_{2}\times _{B}v_{3}\right) =Q_{B}\left( v_{2}\times _{B}\left(
-v_{1}-v_{2}\right) \right) =Q_{B}\left( -v_{2}\times _{B}v_{1}\right)
=Q_{B}\left( v_{1}\times v_{2}\right)
\end{equation*}%
and similarly 
\begin{equation*}
Q_{B}\left( v_{2}\times _{B}v_{3}\right) =Q_{B}\left( v_{3}\times
_{B}v_{1}\right) .
\end{equation*}
From Lagrange's identity 
\begin{equation*}
Q_{B}\left( v_{1}\times _{B}v_{2}\right) =\left( \det B\right) \left[
Q_{B}\left( v_{1}\right) Q_{B}\left( v_{2}\right) -\left( v_{1}\cdot
_{B}v_{2}\right) ^{2}\right]
\end{equation*}%
we use the polarisation formula to obtain%
\begin{eqnarray*}
Q_{B}\left( v_{1}\times _{B}v_{2}\right) &=&\left( \det B\right) \left[
Q_{1}Q_{2}-\left( \frac{Q_{1}+Q_{2}-Q_{3}}{2}\right) ^{2}\right] \\
&=&\frac{\det B}{4}\left[ 4Q_{1}Q_{2}-\left( Q_{1}+Q_{2}-Q_{3}\right) ^{2}%
\right] \\
&=&\frac{\det B}{4}A\left( Q_{1},Q_{2},Q_{3}\right)
\end{eqnarray*}%
and then a simple rearrangement gives one of our desired results; the others
follow by symmetry.
\end{proof}

Over a general field the angle between vectors is not well-defined. So in
rational trigonometry we replace angle with the algebraic quantity called a
spread which is more directly linked to the underlying scalar products. The $%
B$-\textbf{spread} between vectors $v_{1}$ and $v_{2}$ is the number 
\begin{equation*}
s_{B}\left( v_{1},v_{2}\right) \equiv 1-\frac{\left( v_{1}\cdot
_{B}v_{2}\right) ^{2}}{Q_{B}\left( v_{1}\right) Q_{B}\left( v_{2}\right) }.
\end{equation*}%
By Lagrange's identity this can be rewritten as%
\begin{equation*}
s_{B}\left( v_{1},v_{2}\right) =\frac{Q_{B}\left( v_{1}\times
_{B}v_{2}\right) }{\left( \det B\right) Q_{B}\left( v_{1}\right) Q_{B}\left(
v_{2}\right) }.
\end{equation*}%
The $B$-spread is invariant under scalar multiplication of either $v_{1}$ or 
$v_{2}$ (or both). If one or both of $v_{1},v_{2}$ are null vectors then the
spread is undefined.

In what follows, we will consider a vector triangle $\overline{%
v_{1}v_{2}v_{3}}$ with $B$-quadrances%
\begin{equation*}
Q_{1}\equiv Q_{B}\left( v_{1}\right) ,\quad Q_{2}\equiv Q_{B}\left(
v_{2}\right) \quad \mathrm{and\quad }Q_{3}\equiv Q_{B}\left( v_{3}\right)
\end{equation*}%
as well as $B$-spreads%
\begin{equation*}
s_{1}\equiv s_{B}\left( v_{2},v_{3}\right) ,\quad s_{2}\equiv s_{B}\left(
v_{1},v_{3}\right) \quad \mathrm{and\quad }s_{3}\equiv s_{B}\left(
v_{1},v_{2}\right)
\end{equation*}%
and $B$-quadrea%
\begin{equation*}
\mathcal{A}\equiv \mathcal{A}_{B}\left( \overline{v_{1}v_{2}v_{3}}\right) .
\end{equation*}

We now present some results of planar rational trigonometry which connect
these fundamental quantities in three dimensions, with proof. These are
generalisations of the results in \cite[pp. 89-90]{WildDP} to an arbitrary
symmetric bilinear form and to three dimensions.

\begin{theorem}[Cross law]
For a vector triangle $\overline{v_{1}v_{2}v_{3}}$ with $B$-quadrances $%
Q_{1} $, $Q_{2}$ and $Q_{3}$, and corresponding $B$-spreads $s_{1}$, $s_{2}$
and $s_{3}$, we have%
\begin{equation*}
\left( Q_{2}+Q_{3}-Q_{1}\right) ^{2}=4Q_{2}Q_{3}\left( 1-s_{1}\right)
\end{equation*}%
as well as%
\begin{equation*}
\left( Q_{1}+Q_{3}-Q_{2}\right) ^{2}=4Q_{1}Q_{3}\left( 1-s_{2}\right)
\end{equation*}%
and 
\begin{equation*}
\left( Q_{1}+Q_{2}-Q_{3}\right) ^{2}=4Q_{1}Q_{2}\left( 1-s_{3}\right) .
\end{equation*}
\end{theorem}

\begin{proof}
We just prove the last formula, the others follow by symmetry. Rearrange the
polarisation formula to get%
\begin{equation*}
Q_{1}+Q_{2}-Q_{3}=2\left( v_{1}\cdot _{B}v_{2}\right)
\end{equation*}%
and then square both sides to obtain%
\begin{equation*}
\left( Q_{1}+Q_{2}-Q_{3}\right) ^{2}=4\left( v_{1}\cdot _{B}v_{2}\right)
^{2}.
\end{equation*}%
Rearrange the definition of $s_{3}$ to obtain%
\begin{equation*}
\left( v_{1}\cdot _{B}v_{2}\right) ^{2}=Q_{B}\left( v_{1}\right) Q_{B}\left(
v_{2}\right) \left( 1-s_{3}\right) =Q_{1}Q_{2}\left( 1-s_{3}\right) .
\end{equation*}%
Putting these together we get%
\begin{equation*}
\left( Q_{1}+Q_{2}-Q_{3}\right) ^{2}=4Q_{1}Q_{2}\left( 1-s_{3}\right) .
\end{equation*}
\end{proof}

We use the Cross law as a fundamental building block for a number of other
results. For instance, we can express the $B$-quadrea of a triangle in terms
of its $B$-quadrances and $B$-spreads.

\begin{theorem}[Quadrea spread theorem]
For a vector triangle $\overline{v_{1}v_{2}v_{3}}$ with $B$-quadrances $%
Q_{1} $, $Q_{2}$ and $Q_{3}$, corresponding $B$-spreads $s_{1}$, $s_{2}$ and 
$s_{3} $, and $B$-quadrea $\mathcal{A}$, we have 
\begin{equation*}
\mathcal{A}=4Q_{1}Q_{2}s_{3}=4Q_{1}Q_{3}s_{2}=4Q_{2}Q_{3}s_{1}.
\end{equation*}
\end{theorem}

\begin{proof}
Rearrange the equation%
\begin{equation*}
\left( Q_{1}+Q_{2}-Q_{3}\right) ^{2}=4Q_{1}Q_{2}\left( 1-s_{3}\right)
\end{equation*}%
from the Cross law as%
\begin{equation*}
4Q_{1}Q_{2}-\left( Q_{1}+Q_{2}-Q_{3}\right) ^{2}=4Q_{1}Q_{2}s_{3}.
\end{equation*}%
We recognize on the left hand side an asymmetric form of Archimedes'
function, so that%
\begin{equation*}
\left( Q_{1}+Q_{2}+Q_{3}\right) ^{2}-2\left(
Q_{1}^{2}+Q_{1}^{2}+Q_{1}^{2}\right) =4Q_{1}Q_{2}s_{3}
\end{equation*}%
which is%
\begin{equation*}
A\left( Q_{1},Q_{2},Q_{3}\right) =\mathcal{A}=4Q_{1}Q_{2}s_{3}.
\end{equation*}
\end{proof}

We can use the Quadrea spread theorem to determine whether a vector triangle 
$\overline{v_{1}v_{2}v_{3}}$ is \textbf{degenerate}, i.e. when the vectors $%
v_{1},v_{2}$ and $v_{3}$ are collinear.

\begin{theorem}[Triple quad formula]
Let $\overline{v_{1}v_{2}v_{3}}$ be a vector triangle with $B$-quadrances $%
Q_{1}$, $Q_{2}$ and $Q_{3}$. If $\overline{v_{1}v_{2}v_{3}}$ is degenerate
then%
\begin{equation*}
\left( Q_{1}+Q_{2}+Q_{3}\right) ^{2}=2\left(
Q_{1}^{2}+Q_{2}^{2}+Q_{3}^{2}\right) .
\end{equation*}
\end{theorem}

\begin{proof}
If $\overline{v_{1}v_{2}v_{3}}$ is degenerate, then we may suppose that $%
v_{2}\equiv \lambda v_{1}$, for some $\lambda $ in $\mathbb{F}$, so that $%
v_{3}=-\left( 1+\lambda \right) v_{1}$. Thus, by the properties of the $B$%
-quadratic form,%
\begin{equation*}
Q_{2}=\lambda ^{2}Q_{1}\quad \mathrm{and}\quad Q_{3}=\left( 1+\lambda
\right) ^{2}Q_{1}.
\end{equation*}%
So,%
\begin{eqnarray*}
&&\left( Q_{1}+Q_{2}+Q_{3}\right) ^{2}-2\left(
Q_{1}^{2}+Q_{2}^{2}+Q_{3}^{2}\right) \\
&=&\left[ \left( 1+\lambda ^{2}+\left( 1+\lambda \right) ^{2}\right)
^{2}-2\left( 1+\lambda ^{4}+\left( 1+\lambda \right) ^{4}\right) \right]
Q_{1}^{2} \\
&=&\left[ 4\left( \lambda ^{2}+\lambda +1\right) ^{2}-4\left( \lambda
^{2}+\lambda +1\right) ^{2}\right] Q_{1}^{2} \\
&=&0.
\end{eqnarray*}%
The result immediately follows.
\end{proof}

The Cross law also gives the most important result in geometry and
trigonometry: Pythagoras' theorem.

\begin{theorem}[Pythagoras' theorem]
For a vector triangle $\overline{v_{1}v_{2}v_{3}}$ with $B$-quadrances $%
Q_{1} $, $Q_{2}$ and $Q_{3}$, and corresponding $B$-spreads $s_{1}$, $s_{2}$
and $s_{3},$ we have $s_{3}=1$ precisely when 
\begin{equation*}
Q_{1}+Q_{2}=Q_{3}.
\end{equation*}
\end{theorem}

\begin{proof}
The Cross law relation%
\begin{equation*}
\left( Q_{1}+Q_{2}-Q_{3}\right) ^{2}=4Q_{1}Q_{2}\left( 1-s_{3}\right)
\end{equation*}%
together with the assumption that $Q_{1}$ and $Q_{2}$ are non-zero implies
that $s_{3}=1$ precisely when%
\begin{equation*}
Q_{1}+Q_{2}=Q_{3}.
\end{equation*}
\end{proof}

One other use of the Quadrea spread theorem is in determining ratios between 
$B$-spreads and $B$-quadrances, which is a rational analog of the sine law
in classical trigonometry.

\begin{theorem}[Spread law]
For a vector triangle $\overline{v_{1}v_{2}v_{3}}$ with $B$-quadrances $%
Q_{1} $, $Q_{2}$ and $Q_{3}$, and corresponding $B$-spreads $s_{1}$, $s_{2}$
and $s_{3}$, and $B$-quadrea $\mathcal{A}$, we have 
\begin{equation*}
\frac{s_{1}}{Q_{1}}=\frac{s_{2}}{Q_{2}}=\frac{s_{3}}{Q_{3}}=\frac{\mathcal{A}%
}{4Q_{1}Q_{2}Q_{3}}.
\end{equation*}
\end{theorem}

\begin{proof}
Rearrange the Quadrea spread theorem to get%
\begin{equation*}
s_{1}=\frac{\mathcal{A}}{4Q_{2}Q_{3}},\quad s_{2}=\frac{\mathcal{A}}{%
4Q_{1}Q_{3}}\quad \mathrm{and\quad }s_{3}=\frac{\mathcal{A}}{4Q_{1}Q_{2}}.
\end{equation*}%
If the three $B$-spreads are defined, then necessarily all three $B$%
-quadrances are non-zero. So divide $s_{1}$, $s_{2}$ and $s_{3}$ by $Q_{1}$, 
$Q_{2}$ and $Q_{3}$ respectively to get%
\begin{equation*}
\frac{s_{1}}{Q_{1}}=\frac{s_{2}}{Q_{2}}=\frac{s_{3}}{Q_{3}}=\frac{\mathcal{A}%
}{4Q_{1}Q_{2}Q_{3}}
\end{equation*}%
as required.
\end{proof}

Finally we present a result that gives a relationship between the three $B$%
-spreads of a triangle, following the proof in \cite[pp. 89-90]{WildDP}.

\begin{theorem}[Triple spread formula]
For a vector triangle $\overline{v_{1}v_{2}v_{3}}$ with $B$-spreads $s_{1}$, 
$s_{2}$ and $s_{3}$, we have%
\begin{equation*}
\left( s_{1}+s_{2}+s_{3}\right) ^{2}=2\left(
s_{1}^{2}+s_{2}^{2}+s_{3}^{2}\right) +4s_{1}s_{2}s_{3}.
\end{equation*}
\end{theorem}

\begin{proof}
If one of the $B$-spreads is $0$, then from the Spread law they will all be
zero, and likewise with $\mathcal{A}$; thus the formula is immediate in this
case. Otherwise, with the $B$-quadrances $Q_{1}$, $Q_{2}$ and $Q_{3}$ of $%
\overline{v_{1}v_{2}v_{3}}$ as previously defined, the Spread law allows us
to define the non-zero quantity%
\begin{equation*}
D\equiv \frac{4Q_{1}Q_{2}Q_{3}}{\mathcal{A}}
\end{equation*}%
so that%
\begin{equation*}
Q_{1}=Ds_{1},\quad Q_{2}=Ds_{2}\quad \mathrm{and}\text{\quad }Q_{3}=Ds_{3}.
\end{equation*}%
We substitute these into the Cross law and pull out common factors to get%
\begin{equation*}
D^{2}\left( s_{1}+s_{2}-s_{3}\right) ^{2}=4D^{2}s_{1}s_{2}\left(
1-s_{3}\right) .
\end{equation*}%
Now divide by $D^{2}$ and rearrange to obtain%
\begin{equation*}
4s_{1}s_{2}-\left( s_{1}+s_{2}-s_{3}\right) ^{2}=4s_{1}s_{2}s_{3}.
\end{equation*}%
We use the identity discussed earlier in the context of Archimedes' function
to rearrange this to get 
\begin{equation*}
\left( s_{1}+s_{2}+s_{3}\right) ^{2}=2\left(
s_{1}^{2}+s_{2}^{2}+s_{3}^{2}\right) +4s_{1}s_{2}s_{3}.
\end{equation*}
\end{proof}

\section{Rational trigonometry for a projective triangle}

Rational trigonometry has an affine and projective version. The projective
version is typically more algebraically involved. The distinction was first
laid out in \cite{WildUHG1} by framing hyperbolic geometry in a projective
setting, and is generally summarised in \cite{WildAP}. So, the projective
results are the essential formulas for the rational trigonometric approach
to both hyperbolic and spherical, or elliptic, trigonometry. In this paper,
the spherical or elliptic interpretation is primary, as it is most easily
accessible from a Euclidean orientation.

A \textbf{projective vector} $p=\left[ v\right] $ is an expression involving
a non-zero vector $v$ with the convention that 
\begin{equation*}
\left[ v\right] =\left[ \lambda v\right]
\end{equation*}%
for any non-zero number $\lambda .$ If $v=\left( a,b,c\right) $ then we will
write $\left[ v\right] =\left[ a:b:c\right] $ since it is only the
proportion between these three numbers that is important. Clearly a
projective vector can be identified with a one-dimensional subspace of $%
V^{3},$ but it will not be necessary to do so.

A \textbf{projective triangle, }or\textbf{\ tripod, }is a set of three
distinct projective vectors, namely $\left\{ p_{1},p_{2},p_{3}\right\}
=\left\{ \left[ v_{1}\right] ,\left[ v_{2}\right] ,\left[ v_{3}\right]
\right\} $ which we will write as $\overline{p_{1}p_{2}p_{3}}$. Such a
projective triangle is \textbf{degenerate} precisely when $v_{1}$, $v_{2}$
and $v_{3}$ are linearly dependent.

If $p_{1}=\left[ v_{1}\right] $ and $p_{2}=\left[ v_{2}\right] $ are
projective points, then we define the $B$\textbf{-normal} of $p_{1}$ and $%
p_{2}$ to be the projective point%
\begin{equation*}
p_{1}\times _{B}p_{2}\equiv \left[ v_{1}\times _{B}v_{2}\right]
\end{equation*}%
and this is clearly well-defined. We now define the $B$-\textbf{dual} of the
non-degenerate tripod $\overline{p_{1}p_{2}p_{3}}$ to be the tripod $%
\overline{r_{1}r_{2}r_{3}}$, where%
\begin{equation*}
r_{1}\equiv p_{2}\times _{B}p_{3},\quad r_{2}\equiv p_{1}\times
_{B}p_{3}\quad \mathrm{and}\text{\quad }r_{3}\equiv p_{1}\times _{B}p_{2}.
\end{equation*}%
Such a tripod will also be called the $B$-\textbf{dual projective triangle} 
\cite{WildUHG1} of the projective triangle $\overline{p_{1}p_{2}p_{3}}$. The
following result highlights the two-fold symmetry of such a concept.

\begin{theorem}
If the $B$-dual of the non-degenerate tripod $\overline{p_{1}p_{2}p_{3}}$ is 
$\overline{r_{1}r_{2}r_{3}}$, then the $B$-dual of the tripod $\overline{%
r_{1}r_{2}r_{3}}$ is $\overline{p_{1}p_{2}p_{3}}$.
\end{theorem}

\begin{proof}
Let $p_{1}\equiv \left[ v_{1}\right] $, $p_{2}\equiv \left[ v_{2}\right] $
and $p_{3}\equiv \left[ v_{3}\right] $, so that%
\begin{equation*}
r_{1}=p_{2}\times _{B}p_{3}=\left[ v_{2}\times _{B}v_{3}\right] ,\quad
r_{2}=p_{1}\times _{B}p_{3}=\left[ v_{1}\times _{B}v_{3}\right]
\end{equation*}%
and%
\begin{equation*}
r_{3}=p_{1}\times _{B}p_{2}=\left[ v_{1}\times _{B}v_{2}\right] .
\end{equation*}%
Suppose that the $B$-dual of $\overline{r_{1}r_{2}r_{3}}$ is given by $%
\overline{t_{1}t_{2}t_{3}}$, where%
\begin{equation*}
t_{1}\equiv r_{2}\times _{B}r_{3},\quad t_{2}\equiv r_{1}\times
_{B}r_{3}\quad \text{and\quad }t_{3}\equiv r_{1}\times _{B}r_{2}.
\end{equation*}%
By the definition of the $B$-normal, we use Corollary 13 to get 
\begin{eqnarray*}
t_{1} &=&\left[ \left( v_{1}\times _{B}v_{3}\right) \times _{B}\left(
v_{1}\times _{B}v_{2}\right) \right] \\
&=&\left[ \left( \det B\right) \left[ v_{1},v_{3},v_{2}\right] _{B}v_{1}%
\right] .
\end{eqnarray*}%
Since $B$ is non-degenerate and the vectors $v_{1}$, $v_{2}$ and $v_{3}$ are
linearly independent, $\left( \det B\right) \left[ v_{1},v_{3},v_{2}\right]
_{B}\neq 0$ and by the definition of a projective point%
\begin{equation*}
t_{1}=\left[ v_{1}\right] =p_{1}.
\end{equation*}%
By symmetry, $t_{2}=p_{2}$ and $t_{3}=p_{3}$, and hence $\overline{%
t_{1}t_{2}t_{3}}=\overline{p_{1}p_{2}p_{3}}$. Thus, the $B$-dual of $%
\overline{r_{1}r_{2}r_{3}}$ is $\overline{p_{1}p_{2}p_{3}}$.
\end{proof}

The $B$\textbf{-projective quadrance} between two projective vectors $p_{1}=%
\left[ v_{1}\right] $ and $p=\left[ v_{2}\right] $ is%
\begin{equation*}
q_{B}\left( p_{1},p_{2}\right) \equiv s_{B}\left( v_{1},v_{2}\right) .
\end{equation*}

Clearly a $B$-projective quadrance is just the $B$-spread between the
corresponding vectors. So, the following result should not be a surprise.

\begin{theorem}[Projective triple quad formula]
If $\overline{p_{1}p_{2}p_{3}}$ is a degenerate tripod with projective
quadrances 
\begin{equation*}
q_{1}\equiv q_{B}\left( p_{2},p_{3}\right) ,\quad q_{2}\equiv q_{B}\left(
p_{1},p_{3}\right) \quad \mathrm{and}\quad q_{3}\equiv q_{B}\left(
p_{1},p_{2}\right) ,
\end{equation*}%
then 
\begin{equation*}
\left( q_{1}+q_{2}+q_{3}\right) ^{2}=2\left(
q_{1}^{2}+q_{2}^{2}+q_{3}^{2}\right) +4q_{1}q_{2}q_{3}.
\end{equation*}
\end{theorem}

The Projective triple quad formula is analogous and parallel to the Triple
spread formula in affine rational trigonometry, due to this fact.

Given a projective triangle $\overline{p_{1}p_{2}p_{3}}$ it will have three
projective quadrances $q_{1},q_{2}$ and $q_{3}$. We now introduce the 
\textbf{projective spreads} $S_{1},S_{2}$ and $S_{3}$ of the projective
triangle $\overline{p_{1}p_{2}p_{3}}$ to be the projective quadrances of its 
$B$-dual $\overline{r_{1}r_{2}r_{3}}$, that is%
\begin{equation*}
S_{1}\equiv q_{B}\left( r_{2},r_{3}\right) ,\quad S_{2}\equiv q_{B}\left(
r_{1},r_{3}\right) \quad \mathrm{and}\quad S_{3}\equiv q_{B}\left(
r_{1},r_{2}\right) .
\end{equation*}

We now proceed to present results in projective rational trigonometry, which
draw on the results from \cite{WildAP}, but will be framed in the
three-dimensional framework using $B$-vector products and a general
symmetric bilinear form.

\begin{theorem}[Projective spread law]
Given a tripod $\overline{p_{1}p_{2}p_{3}}$ with $B$-projective quadrances $%
q_{1}$, $q_{2}$ and $q_{3}$, and $B$-projective spreads $S_{1}$, $S_{2}$ and 
$S_{3}$, we have%
\begin{equation*}
\frac{S_{1}}{q_{1}}=\frac{S_{2}}{q_{2}}=\frac{S_{3}}{q_{3}}.
\end{equation*}
\end{theorem}

\begin{proof}
Let $p_{1}\equiv \left[ v_{1}\right] $, $p_{2}\equiv \left[ v_{2}\right] $
and $p_{3}\equiv \left[ v_{3}\right] $ be the points of $\overline{%
p_{1}p_{2}p_{3}}$. Also consider its $B$-dual $\overline{r_{1}r_{2}r_{3}}$,
where%
\begin{equation*}
r_{1}\equiv \left[ v_{2}\times _{B}v_{3}\right] ,\quad r_{2}\equiv \left[
v_{1}\times _{B}v_{3}\right] \quad \mathrm{and\quad }r_{3}\equiv \left[
v_{1}\times _{B}v_{2}\right] .
\end{equation*}%
By the definition of the $B$-projective quadrance,%
\begin{equation*}
q_{1}=\frac{Q_{B}\left( v_{2}\times _{B}v_{3}\right) }{\left( \det B\right)
Q_{B}\left( v_{2}\right) Q_{B}\left( v_{3}\right) }
\end{equation*}%
and similarly%
\begin{equation*}
q_{2}=\frac{Q_{B}\left( v_{1}\times _{B}v_{3}\right) }{\left( \det B\right)
Q_{B}\left( v_{1}\right) Q_{B}\left( v_{3}\right) }\quad \mathrm{and}\quad
q_{3}=\frac{Q_{B}\left( v_{1}\times _{B}v_{2}\right) }{\left( \det B\right)
Q_{B}\left( v_{1}\right) Q_{B}\left( v_{2}\right) }.
\end{equation*}%
By the definition of the $B$-projective spread and Corollary 13, 
\begin{equation*}
S_{1}=\frac{Q_{B}\left( \left\langle v_{1},v_{2};v_{1},v_{3}\right\rangle
_{B}\right) }{\left( \det B\right) Q_{B}\left( v_{1}\times _{B}v_{2}\right)
Q_{B}\left( v_{1}\times _{B}v_{3}\right) }=\frac{\left( \det B\right) \left[
v_{1},v_{2},v_{3}\right] _{B}^{2}Q_{B}\left( v_{1}\right) }{Q_{B}\left(
v_{1}\times _{B}v_{2}\right) Q_{B}\left( v_{1}\times _{B}v_{3}\right) }
\end{equation*}%
and similarly%
\begin{equation*}
S_{2}=\frac{\left( \det B\right) \left[ v_{1},v_{2},v_{3}\right]
_{B}^{2}Q_{B}\left( v_{2}\right) }{Q_{B}\left( v_{1}\times _{B}v_{2}\right)
Q_{B}\left( v_{2}\times _{B}v_{3}\right) }\quad \mathrm{and\quad }S_{3}=%
\frac{\left( \det B\right) \left[ v_{1},v_{2},v_{3}\right]
_{B}^{2}Q_{B}\left( v_{3}\right) }{Q_{B}\left( v_{1}\times _{B}v_{3}\right)
Q_{B}\left( v_{2}\times _{B}v_{3}\right) }.
\end{equation*}%
So,%
\begin{equation*}
\frac{S_{1}}{q_{1}}=\frac{\left( \det B\right) ^{2}\left[ v_{1},v_{2},v_{3}%
\right] _{B}^{2}Q_{B}\left( v_{1}\right) Q_{B}\left( v_{2}\right)
Q_{B}\left( v_{3}\right) }{Q_{B}\left( v_{1}\times _{B}v_{2}\right)
Q_{B}\left( v_{1}\times _{B}v_{3}\right) Q_{B}\left( v_{2}\times
_{B}v_{3}\right) }=\frac{S_{2}}{q_{2}}=\frac{S_{3}}{q_{3}},
\end{equation*}%
as required.
\end{proof}

If we balance each side of the result of the Projective spread law to its
lowest common denominator, multiplying through by the denominator motivates
us to define the quantity%
\begin{equation*}
S_{1}q_{2}q_{3}=S_{2}q_{1}q_{3}=S_{3}q_{1}q_{2}\equiv a_{B}\left( \overline{%
p_{1}p_{2}p_{3}}\right) \equiv a_{B},
\end{equation*}%
which will be called the $B$-\textbf{projective quadrea} of the tripod $%
\overline{p_{1}p_{2}p_{3}}$.

There is a relationship between the $B$-projective quadrea and the
projective quadrances discovered in \cite{WildAP}, which is central to our
study of projective rational trigonometry. We extend this result to an
arbitrary symmetric bilinear form, using a quite different argument.

\begin{theorem}[Projective cross law]
Given a tripod $\overline{p_{1}p_{2}p_{3}}$ with $B$-projective quadrances $%
q_{1}$, $q_{2}$ and $q_{3}$, $B$-projective spreads $S_{1}$, $S_{2}$ and $%
S_{3}$, and $B$-projective quadrea $a_{B}$,%
\begin{equation*}
\left( a_{B}-q_{1}-q_{2}-q_{3}+2\right) ^{2}=4\left( 1-q_{1}\right) \left(
1-q_{2}\right) \left( 1-q_{3}\right) .
\end{equation*}
\end{theorem}

\begin{proof}
Let $p_{1}\equiv \left[ v_{1}\right] $, $p_{2}\equiv \left[ v_{2}\right] $
and $p_{3}\equiv \left[ v_{3}\right] $ be the points of $\overline{%
p_{1}p_{2}p_{3}}$ and $\overline{r_{1}r_{2}r_{3}}$ be the $B$-dual of $%
\overline{p_{1}p_{2}p_{3}}$, so that from the proof of the Projective spread
law the $B$-projective quadrances and spreads are%
\begin{equation*}
q_{1}=\frac{Q_{B}\left( v_{2}\times _{B}v_{3}\right) }{\left( \det B\right)
Q_{B}\left( v_{2}\right) Q_{B}\left( v_{3}\right) },\quad q_{2}=\frac{%
Q_{B}\left( v_{1}\times _{B}v_{3}\right) }{\left( \det B\right) Q_{B}\left(
v_{1}\right) Q_{B}\left( v_{3}\right) },\quad q_{3}=\frac{Q_{B}\left(
v_{1}\times _{B}v_{2}\right) }{\left( \det B\right) Q_{B}\left( v_{1}\right)
Q_{B}\left( v_{2}\right) },
\end{equation*}%
\begin{equation*}
S_{1}=\frac{\left( \det B\right) \left[ v_{1},v_{2},v_{3}\right]
_{B}^{2}Q_{B}\left( v_{1}\right) }{Q_{B}\left( v_{1}\times _{B}v_{2}\right)
Q_{B}\left( v_{1}\times _{B}v_{3}\right) },\quad S_{2}=\frac{\left( \det
B\right) \left[ v_{1},v_{2},v_{3}\right] _{B}^{2}Q_{B}\left( v_{2}\right) }{%
Q_{B}\left( v_{1}\times _{B}v_{2}\right) Q_{B}\left( v_{2}\times
_{B}v_{3}\right) }
\end{equation*}%
and%
\begin{equation*}
S_{3}=\frac{\left( \det B\right) \left[ v_{1},v_{2},v_{3}\right]
_{B}^{2}Q_{B}\left( v_{3}\right) }{Q_{B}\left( v_{1}\times _{B}v_{3}\right)
Q_{B}\left( v_{2}\times _{B}v_{3}\right) }.
\end{equation*}%
Furthermore,%
\begin{eqnarray*}
a_{B} &=&S_{1}q_{2}q_{3}=S_{2}q_{1}q_{3}=S_{3}q_{1}q_{2} \\
&=&\frac{\left[ v_{1},v_{2},v_{3}\right] _{B}^{2}}{\left( \det B\right)
Q_{B}\left( v_{1}\right) Q_{B}\left( v_{2}\right) Q_{B}\left( v_{3}\right) }.
\end{eqnarray*}%
Noting that%
\begin{equation*}
1-q_{1}=\frac{\left( v_{2}\cdot _{B}v_{3}\right) ^{2}}{Q_{B}\left(
v_{2}\right) Q_{B}\left( v_{3}\right) },\quad 1-q_{2}=\frac{\left(
v_{1}\cdot _{B}v_{3}\right) ^{2}}{Q_{B}\left( v_{1}\right) Q_{B}\left(
v_{3}\right) }\quad \mathrm{and}\quad 1-q_{3}=\frac{\left( v_{1}\cdot
_{B}v_{2}\right) ^{2}}{Q_{B}\left( v_{1}\right) Q_{B}\left( v_{2}\right) }
\end{equation*}%
we see from the Scalar triple product theorem that 
\begin{eqnarray*}
&&\left( a_{B}-q_{1}-q_{2}-q_{3}+2\right) ^{2} \\
&=&\left( \frac{\left( \det M\right) ^{2}\det B}{Q_{B}\left( v_{1}\right)
Q_{B}\left( v_{2}\right) Q_{B}\left( v_{3}\right) }+\left( 1-q_{1}\right)
+\left( 1-q_{2}\right) +\left( 1-q_{3}\right) -1\right) ^{2} \\
&=&\left( \frac{\left( \det M\right) ^{2}\det B}{Q_{B}\left( v_{1}\right)
Q_{B}\left( v_{2}\right) Q_{B}\left( v_{3}\right) }+\frac{\left( v_{2}\cdot
_{B}v_{3}\right) ^{2}}{Q_{B}\left( v_{2}\right) Q_{B}\left( v_{3}\right) }+%
\frac{\left( v_{1}\cdot _{B}v_{3}\right) ^{2}}{Q_{B}\left( v_{1}\right)
Q_{B}\left( v_{3}\right) }+\frac{\left( v_{1}\cdot _{B}v_{2}\right) ^{2}}{%
Q_{B}\left( v_{1}\right) Q_{B}\left( v_{2}\right) }-1\right) ^{2}
\end{eqnarray*}%
where%
\begin{equation*}
M\equiv 
\begin{pmatrix}
- & v_{1} & - \\ 
- & v_{2} & - \\ 
- & v_{3} & -%
\end{pmatrix}%
.
\end{equation*}%
Given that%
\begin{eqnarray*}
\left( \det M\right) ^{2}\det B &=&\det \left( MBM^{T}\right) \\
&=&\det 
\begin{pmatrix}
Q_{B}\left( v_{1}\right) & v_{1}\cdot _{B}v_{2} & v_{1}\cdot _{B}v_{3} \\ 
v_{1}\cdot _{B}v_{2} & Q_{B}\left( v_{2}\right) & v_{2}\cdot _{B}v_{3} \\ 
v_{1}\cdot _{B}v_{3} & v_{2}\cdot _{B}v_{3} & Q_{B}\left( v_{3}\right)%
\end{pmatrix}
\\
&=&Q_{B}\left( v_{1}\right) Q_{B}\left( v_{2}\right) Q_{B}\left(
v_{3}\right) +2\left( v_{1}\cdot _{B}v_{2}\right) \left( v_{1}\cdot
_{B}v_{3}\right) \left( v_{2}\cdot _{B}v_{3}\right) \\
&&-\left( v_{1}\cdot _{B}v_{2}\right) ^{2}Q_{B}\left( v_{3}\right) -\left(
v_{1}\cdot _{B}v_{3}\right) ^{2}Q_{B}\left( v_{2}\right) -\left( v_{2}\cdot
_{B}v_{3}\right) ^{2}Q_{B}\left( v_{1}\right)
\end{eqnarray*}%
we obtain%
\begin{eqnarray*}
\left( a_{B}-q_{1}-q_{2}-q_{3}+2\right) ^{2} &=&\left( \frac{2\left(
v_{1}\cdot _{B}v_{2}\right) \left( v_{1}\cdot _{B}v_{3}\right) \left(
v_{2}\cdot _{B}v_{3}\right) }{Q_{B}\left( v_{1}\right) Q_{B}\left(
v_{2}\right) Q_{B}\left( v_{3}\right) }\right) ^{2} \\
&=&4\frac{\left( v_{1}\cdot _{B}v_{2}\right) ^{2}}{Q_{B}\left( v_{1}\right)
Q_{B}\left( v_{2}\right) }\frac{\left( v_{1}\cdot _{B}v_{3}\right) ^{2}}{%
Q_{B}\left( v_{1}\right) Q_{B}\left( v_{3}\right) }\frac{\left( v_{2}\cdot
_{B}v_{3}\right) ^{2}}{Q_{B}\left( v_{2}\right) Q_{B}\left( v_{3}\right) } \\
&=&4\left( 1-q_{1}\right) \left( 1-q_{2}\right) \left( 1-q_{3}\right)
\end{eqnarray*}%
as required.
\end{proof}

The Projective cross law can also be expressed in various asymmetric forms.

\begin{corollary}
Given a tripod $\overline{p_{1}p_{2}p_{3}}$ with $B$-projective quadrances $%
q_{1}$, $q_{2}$ and $q_{3}$, and $B$-projective spreads $S_{1}$, $S_{2}$ and 
$S_{3}$, the Projective cross law can be rewritten as either%
\begin{equation*}
\left( S_{1}q_{2}q_{3}+q_{1}-q_{2}-q_{3}\right) ^{2}=4q_{2}q_{3}\left(
1-q_{1}\right) \left( 1-S_{1}\right) ,
\end{equation*}%
\begin{equation*}
\left( S_{2}q_{1}q_{3}-q_{1}+q_{2}-q_{3}\right) ^{2}=4q_{1}q_{3}\left(
1-q_{2}\right) \left( 1-S_{2}\right)
\end{equation*}%
or%
\begin{equation*}
\left( S_{3}q_{1}q_{2}-q_{1}-q_{2}+q_{3}\right) ^{2}=4q_{1}q_{2}\left(
1-q_{3}\right) \left( 1-S_{3}\right) .
\end{equation*}
\end{corollary}

\begin{proof}
Let $C_{1}\equiv 1-S_{1}$. Substitute $a_{B}=S_{1}q_{2}q_{3}=\left(
1-C_{1}\right) q_{2}q_{3}$ into the Projective cross law to get%
\begin{equation*}
\left( \left( 1-C_{1}\right) q_{2}q_{3}-q_{1}-q_{2}-q_{3}+2\right)
^{2}-4\left( 1-q_{1}\right) \left( 1-q_{2}\right) \left( 1-q_{3}\right) =0.
\end{equation*}%
Expand the left-hand side and simplify the result as a polynomial in $C_{1}$
to obtain%
\begin{equation*}
\allowbreak \left( q_{2}^{2}q_{3}^{2}\right) C_{1}^{2}+2q_{2}q_{3}\left(
q_{1}+q_{2}+q_{3}-q_{2}q_{3}-2\right) \allowbreak C_{1}+\left(
q_{2}q_{3}+q_{1}-q_{2}-q_{3}\right) ^{2}=0.
\end{equation*}%
Now insert the term $-4q_{12}q_{13}q_{23}C_{1}$ into the last equation and
balance the equation as required. Further rearrange to get%
\begin{eqnarray*}
&&\left( q_{2}^{2}q_{3}^{2}\right) C_{1}^{2}+2q_{2}q_{3}\left(
q_{2}+q_{3}-q_{1}-q_{2}q_{3}\right) \allowbreak C_{1}+\left(
q_{2}q_{3}+q_{1}-q_{2}-q_{3}\right) ^{2} \\
&=&4q_{2}q_{3}C_{1}-4q_{1}q_{2}q_{3}C_{1}=4q_{2}q_{3}\left( 1-q_{1}\right)
C_{1}.
\end{eqnarray*}%
As the left-hand side is a perfect square, factorise this to get%
\begin{equation*}
\left( q_{2}+q_{3}-q_{1}-q_{2}q_{3}+q_{2}q_{3}C_{1}\right)
^{2}=4q_{2}q_{3}\left( 1-q_{1}\right) C_{1}.
\end{equation*}%
Replace $C_{1}$ with $1-S_{1}$ and simplify to obtain%
\begin{eqnarray*}
&&\left( q_{2}+q_{3}-q_{1}-q_{2}q_{3}+\left( 1-S_{1}\right)
q_{2}q_{3}\right) ^{2} \\
&=&\left( q_{2}q_{3}S_{1}+q_{1}-q_{2}-q_{3}\right) ^{2}=4q_{2}q_{3}\left(
1-q_{1}\right) \left( 1-S_{1}\right) .
\end{eqnarray*}%
The other results follow by symmetry.
\end{proof}

Note that the $B$-projective quadrea $a_{B}$ also features in the
reformulation, and can replace the quantity $S_{1}q_{2}q_{3}$ (as well as
its symmetrical reformulations). Also, the Projective triple quad formula
from earlier follows directly from the Projective cross law, so it can be
proven in this way; this will be omitted from the paper.

In addition to the $B$-projective quadrea, we can also discuss the dual
analog of it. This quantity is called the $B$-\textbf{quadreal }\cite%
{WildUHG1} and is defined by%
\begin{equation*}
l_{B}\equiv l_{B}\left( \overline{p_{1}p_{2}p_{3}}\right) \equiv
q_{1}S_{2}S_{3}=q_{2}S_{1}S_{3}=q_{3}S_{1}S_{2}.
\end{equation*}%
We can also say from Theorem 23 that $l_{B}$ is the $B$-projective quadrea
of the $B$-dual tripod $\overline{r_{1}r_{2}r_{3}}$ of $\overline{%
p_{1}p_{2}p_{3}}$ and $a_{B}$ is the $B$-quadreal of $\overline{%
r_{1}r_{2}r_{3}}$. The following extends the result in \cite{WildUHG1} for $%
B $-quadratic forms.

\begin{corollary}
For a tripod $\overline{p_{1}p_{2}p_{3}}$ with $B$-projective quadrances $%
q_{1}$, $q_{2}$ and $q_{3}$, $B$-projective spreads $S_{1}$, $S_{2}$ and $%
S_{3}$, $B$-projective quadrea $a_{B}$ and $B$-quadreal $l_{B}$,%
\begin{equation*}
a_{B}l_{B}=q_{1}q_{2}q_{3}S_{1}S_{2}S_{3}.
\end{equation*}
\end{corollary}

\begin{proof}
Given%
\begin{equation*}
a_{B}=S_{1}q_{2}q_{3}=S_{2}q_{1}q_{3}=S_{3}q_{1}q_{2}
\end{equation*}%
and%
\begin{equation*}
l_{B}=q_{1}S_{2}S_{3}=q_{2}S_{1}S_{3}=q_{3}S_{1}S_{2},
\end{equation*}%
we get%
\begin{eqnarray*}
a_{B}l_{B} &=&\left( S_{1}q_{2}q_{3}\right) \left( q_{1}S_{2}S_{3}\right)
=\left( S_{2}q_{1}q_{3}\right) \left( q_{2}S_{1}S_{3}\right) \\
&=&\left( S_{3}q_{1}q_{2}\right) \left( q_{3}S_{1}S_{2}\right)
=q_{12}q_{13}q_{23}S_{1}S_{2}S_{3},
\end{eqnarray*}%
as required.
\end{proof}

We now present a projective version of Pythagoras' theorem. This is an
extension of the result in \cite{WildAP} and \cite{WildUHG1} to arbitrary
symmetric bilinear forms.

\begin{theorem}[Projective Pythagoras' theorem]
Take a tripod $\overline{p_{1}p_{2}p_{3}}$ with $B$-projective quadrances $%
q_{1}$, $q_{2}$ and $q_{3}$, and $B$-projective spreads $S_{1}$, $S_{2}$ and 
$S_{3}$. If $S_{1}=1$, then%
\begin{equation*}
q_{1}=q_{2}+q_{3}-q_{2}q_{3}.
\end{equation*}
\end{theorem}

\begin{proof}
Substitute $S_{1}=1$ into the Projective cross law%
\begin{equation*}
\left( S_{1}q_{2}q_{3}-q_{1}-q_{2}-q_{3}+2\right) ^{2}=4\left(
1-q_{1}\right) \left( 1-q_{2}\right) \left( 1-q_{3}\right)
\end{equation*}
and rearrange the result to get%
\begin{equation*}
\left( q_{2}q_{3}-q_{1}-q_{2}-q_{3}+2\right) ^{2}-4\left( 1-q_{1}\right)
\left( 1-q_{2}\right) \left( 1-q_{3}\right) =0.
\end{equation*}%
The left-hand side is factored into%
\begin{equation*}
\left( q_{2}+q_{3}-q_{1}-q_{2}q_{3}\right) ^{2}=0,
\end{equation*}%
so that solving for $q_{1}$ gives%
\begin{equation*}
q_{1}=q_{2}+q_{3}-q_{2}q_{3}.
\end{equation*}
\end{proof}

Note the term $-q_{2}q_{3}$ involved in the Projective Pythagoras' theorem;
this is not present in Pythagoras' theorem in affine rational trigonometry.
The Projective Pythagoras' theorem can be restated \cite{WildUHG1} as%
\begin{eqnarray*}
1-q_{1} &=&1-q_{2}-q_{3}+q_{2}q_{3} \\
&=&\left( 1-q_{2}\right) \left( 1-q_{3}\right) .
\end{eqnarray*}%
As for the converse of the Projective Pythagoras' theorem, start with the
asymmetric form of the Projective cross law%
\begin{equation*}
\left( S_{1}q_{2}q_{3}+q_{1}-q_{2}-q_{3}\right) ^{2}=4q_{2}q_{3}\left(
1-q_{1}\right) \left( 1-S_{1}\right) .
\end{equation*}%
If $q_{1}=q_{2}+q_{3}-q_{2}q_{3}$ then%
\begin{equation*}
\left( q_{2}q_{3}\left( 1-S_{1}\right) \right) ^{2}=4q_{2}q_{3}\left(
1-q_{2}\right) \left( 1-q_{3}\right) \left( 1-S_{1}\right) ,
\end{equation*}%
which can also be rearranged and factorised as%
\begin{equation*}
q_{2}q_{3}\left( 1-S_{1}\right) \left(
4q_{2}+4q_{3}-3q_{2}q_{3}-S_{1}q_{2}q_{3}-4\right) =0.
\end{equation*}%
Here, we see that $S_{1}=1$ is not the only solution; we can also have the
solution%
\begin{equation*}
S_{1}=\frac{4\left( q_{2}+q_{3}-1\right) }{q_{2}q_{3}}-3.
\end{equation*}%
So, the converse of the Projective Pythagoras' theorem may not necessarily
hold. It is of independent interest to deduce the meaning of the latter
solution for $S_{1}$.

\section{The methane molecule for chemists}

To illustrate the practical aspect and attractive values that this
technology gives, we apply the results of this paper to the methane molecule 
$CH_{4}$ consisting of four hydrogen atoms arranged in the form of a regular
tetrahedron, and a central carbon atom. The geometry of this configuration
is well-known to chemists, at least using the classical measurements;
however with rational trigonometry a new picture emerges which illustrates
the advantages of thinking algebraically.

Note that we do not assume a particular field here; to make things precise
mathematically we would need an appropriate quadratic extension of the
rationals to fix the vectors in an appropriate vector space, but we do work
over the familiar Euclidean geometry, so we remove the $B$ prefix from the
subsequent discussion. Once we have built the regular tetrahedron, all the
measurements are rational expressions in the common quadrance of the six
sides, which we will denote by $Q.$ Then the faces are equilateral triangles
with quadreas%
\begin{equation*}
\mathcal{A}=A\left( Q,Q,Q\right) =\left( 3Q\right) ^{2}-2(3Q^{2})=3Q^{2}
\end{equation*}%
which is $16$ times the square of the classical area. The spreads $s$ in any
such equilateral triangle satisfy the Cross law 
\begin{equation*}
\left( Q+Q-Q\right) ^{2}=4Q\times Q\times \left( 1-s\right)
\end{equation*}%
so that 
\begin{equation*}
s=\frac{3}{4}
\end{equation*}%
which is the (rational) analog of approximately $1.\,\allowbreak 047\,20$ in
the radian system, or exactly $60^{\circ }$ in the much more ancient
Babylonian system.

This will also be the projective quadrance $q$ of any two sides meeting at a
common vertex, so $q=3/4$. Three such concurrent sides gives an equilateral
projective triangle, and the projective formulas we have developed apply
also to the (equal) projective spreads $S$ of this projective triangle. In
particular the Projective cross law in terms of the projective quadrea $a$
gives%
\begin{equation*}
\left( a-3\times \frac{3}{4}+2\right) ^{2}=4\left( 1-\frac{3}{4}\right) ^{3}
\end{equation*}%
or 
\begin{equation*}
\left( a-\frac{1}{4}\right) ^{2}=\frac{1}{16}
\end{equation*}%
and since $a\neq 0$ we get 
\begin{equation*}
a=\frac{1}{2}.
\end{equation*}%
This quantity can be considered as the solid spread formed by the three
lines meeting at a vertex, which is a rational analog of the solid angle of
spherical trigonometry. But from the definition of the projective quadrea,
and the symmetry of the Projective spread law, we deduce that $a=Sq^{2}$ so
that%
\begin{equation*}
S=\frac{8}{9}.
\end{equation*}%
Geometrically this is the spread between two faces of the tetrahedron, which
is the (rational) analog of approximately $\arcsin \left( \sqrt{8/9}\right)
\simeq \allowbreak 1.\,\allowbreak 230\,96$ in the radian system, or
approximately $70.\,\allowbreak 528\,8^{\circ }$ in the Babylonian system.
The supplement of this angle, which is approximately\allowbreak\ $%
109.\,\allowbreak 471^{\circ },$ corresponds to the same spread of $8/9,$
which geometrically is formed by the central lines from the carbon atom to
any two hydrogen atom. So we see that the rational trigonometry of this
paper allows more natural, rational and exact expressions for an important
measurement in chemistry, with the calculations in the realm of high school
algebra, without use of a calculator. This is a powerful indicator that for
more complicated calculations, we can expect a significant speed up of
processing by adopting the language and concepts of rational trigonometry.

In our follow up paper we will be investigating the rich trigonometry of a
general tetrahedron, for which this is just a very simple example.

\section{Affine and projective relativistic trigonometry in relativistic
geometries}

To illustrate both affine and projective formulas in a less symmetrical
situation, we shift to a relativistic-three dimensional geometry, and
consider two examples of triangles: one of a vector triangle in $\mathbb{V}%
^{3}$ and the other of a projective triangle in $\mathbb{P}^{2}$, both over
the rational number field. Here the metric structure is given by the \textbf{%
Minkowski scalar product} \cite{Minkowski}\ on $\mathbb{V}^{3}$ defined by%
\begin{equation*}
\left( x_{1},y_{1},z_{1}\right) \cdot _{B}\left( x_{2},y_{2},z_{2}\right)
\equiv x_{1}x_{2}+y_{1}y_{2}-z_{1}z_{2}.
\end{equation*}%
The matrix%
\begin{equation*}
B\equiv 
\begin{pmatrix}
1 & 0 & 0 \\ 
0 & 1 & 0 \\ 
0 & 0 & -1%
\end{pmatrix}%
\end{equation*}%
which represents this symmetric bilinear form, is often called the \textbf{%
relativistic bilinear form}.

\subsection*{Affine relativistic example}

In the first example, consider the vector triangle $\overline{v_{1}v_{2}v_{3}%
}$ with%
\begin{equation*}
v_{1}\equiv \left( -1,3,-2\right) ,\quad v_{2}\equiv \left( 2,-5,4\right) 
\mathrm{\quad and\quad }v_{3}\equiv -v_{1}-v_{2}=\left( -1,2,-2\right) .
\end{equation*}%
The $B$-quadrances of $\overline{v_{1}v_{2}v_{3}}$ are%
\begin{equation*}
Q_{1}=\left( -1,3,-2\right) 
\begin{pmatrix}
1 & 0 & 0 \\ 
0 & 1 & 0 \\ 
0 & 0 & -1%
\end{pmatrix}%
\left( -1,3,-2\right) ^{T}=6,
\end{equation*}%
and similarly%
\begin{equation*}
Q_{2}=13\quad \mathrm{and}\quad Q_{3}=1.
\end{equation*}%
The $B$-quadrea of $\overline{A_{0}A_{1}A_{2}}$ is then%
\begin{equation*}
\mathcal{A}=\left( 6+13+1\right) ^{2}-2\left( 6^{2}+13^{2}+1^{2}\right) =-12,
\end{equation*}%
and hence, by the Quadrea spread theorem, the $B$-spreads are%
\begin{equation*}
s_{1}=\frac{-12}{4\times 13\times 1}=-\frac{3}{13},\quad s_{2}=\frac{-12}{%
4\times 6\times 1}=-\frac{1}{2}\quad \mathrm{and}\text{\quad }s_{3}=\frac{-12%
}{4\times 6\times 13}=-\frac{1}{26}.
\end{equation*}%
To verify our calculations, we observe that%
\begin{equation*}
\frac{s_{1}}{Q_{1}}=\frac{-3}{13\times 6}=-\frac{1}{26},\quad \frac{s_{2}}{%
Q_{2}}=-\frac{1}{13\times 2}=-\frac{1}{26}\quad \mathrm{and}\text{\quad }%
\frac{s_{3}}{Q_{3}}=-\frac{1}{26\times 1}=-\frac{1}{26}.
\end{equation*}%
As each ratio is equal, the Spread law holds for $\overline{A_{0}A_{1}A_{2}}$%
. Furthermore,%
\begin{eqnarray*}
&&\left( s_{0}+s_{1}+s_{2}\right) ^{2}-2\left(
s_{0}^{2}+s_{1}^{2}+s_{2}^{2}\right) \\
&=&\left( -\frac{3}{13}-\frac{1}{2}-\frac{1}{26}\right) ^{2}-2\left( \left( -%
\frac{3}{13}\right) ^{2}+\left( -\frac{1}{2}\right) ^{2}+\left( -\frac{1}{26}%
\right) ^{2}\right) \\
&=&-\frac{3}{169}
\end{eqnarray*}%
and 
\begin{equation*}
4s_{0}s_{1}s_{2}=4\left( -\frac{3}{13}\right) \left( -\frac{1}{2}\right)
\left( -\frac{1}{26}\right) =-\frac{3}{169}.
\end{equation*}%
Because of the equality of these two identities, the Triple spread formula
thus holds.

\subsection*{Projective relativistic example}

Now for a second projective example, let%
\begin{equation*}
v_{1}\equiv \left( 2,-1,3\right) ,\quad v_{2}\equiv \left( -2,5,0\right)
\quad \mathrm{and}\text{\quad }v_{3}\equiv \left( 3,0,4\right)
\end{equation*}%
be three vectors in $\mathbb{V}^{3}$, so that we may define $\overline{%
p_{1}p_{2}p_{3}}$ to be a projective triangle in $\mathbb{P}^{2}$ with
projective points%
\begin{equation*}
p_{1}\equiv \left[ v_{1}\right] ,\quad p_{2}\equiv \left[ v_{2}\right] \quad 
\mathrm{and}\quad p_{3}\equiv \left[ v_{3}\right] .
\end{equation*}%
The $B$-projective quadrances of $\overline{p_{1}p_{2}p_{3}}$ are%
\begin{equation*}
q_{1}=1-\frac{\left( -6\right) ^{2}}{29\times -7}=\frac{239}{203}
\end{equation*}%
$\allowbreak $and similarly%
\begin{equation*}
q_{2}=-\frac{2}{7}\quad \mathrm{and}\text{\quad }q_{3}=\frac{197}{116}.
\end{equation*}%
Let $\overline{r_{1}r_{2}r_{3}}$ be the $B$-dual of $\overline{%
p_{1}p_{2}p_{3}}$, so that%
\begin{equation*}
r_{1}=\left[ \left( 20,8,15\right) \right] ,\quad r_{2}=\left[ \left(
4,-1,3\right) \right] \quad \mathrm{and}\quad r_{3}=\left[ \left(
15,6,8\right) \right] .
\end{equation*}%
Then, the $B$-projective spreads of $\overline{p_{1}p_{2}p_{3}}$ are 
\begin{equation*}
S_{1}=1-\frac{30^{2}}{197\times 8}=\allowbreak \frac{169}{394},
\end{equation*}%
and similarly%
\begin{equation*}
S_{2}=-\frac{4901}{47\,083}\quad \mathrm{and}\text{\quad }S_{3}=\allowbreak 
\frac{1183}{1912}.
\end{equation*}%
To verify our calculations, we observe that%
\begin{equation*}
\frac{S_{1}}{q_{1}}=\allowbreak \frac{169}{394}\div \frac{239}{203}=\frac{%
34\,307}{94\,166},
\end{equation*}%
and similarly%
\begin{equation*}
\frac{S_{2}}{q_{2}}=\left( -\frac{4901}{47\,083}\right) \div \left( -\frac{2%
}{7}\right) =\frac{34\,307}{94\,166}\quad \text{and}\quad \frac{S_{3}}{q_{3}}%
=\frac{1183}{1912}\div \frac{197}{116}=\frac{34\,307}{94\,166}.
\end{equation*}%
Since they are all equal, the Projective spread law holds. Furthermore, with
the $B$-projective quadrea of $\overline{p_{1}p_{2}p_{3}}$ given as%
\begin{equation*}
a_{B}=\frac{197}{116}\times \left( -\frac{2}{7}\right) \times \allowbreak 
\frac{169}{394}=-\frac{169}{812},
\end{equation*}%
we observe that 
\begin{eqnarray*}
&&\left( a_{B}-q_{12}-q_{13}-q_{23}+2\right) ^{2} \\
&=&\left( -\frac{169}{812}-\frac{197}{116}+\frac{2}{7}-\frac{239}{203}%
+2\right) ^{2} \\
&=&\frac{26\,244}{41\,209}
\end{eqnarray*}%
and%
\begin{eqnarray*}
&&4\left( 1-q_{12}\right) \left( 1-q_{13}\right) \left( 1-q_{23}\right) \\
&=&4\left( 1-\frac{197}{116}\right) \left( 1+\frac{2}{7}\right) \left( 1-%
\frac{239}{203}\right) \\
&=&\frac{26\,244}{41\,209}.
\end{eqnarray*}%
As we have equality, the Projective cross law is verified. Note that for the 
$B$-dual tripod $\overline{r_{1}r_{2}r_{3}}$ of $\overline{p_{1}p_{2}p_{3}}$%
, its $B$-projective quadrances are the $B$-projective spreads of $\overline{%
p_{1}p_{2}p_{3}}$ and its $B$-projective spreads are the $B$-projective
quadrances of $\overline{p_{1}p_{2}p_{3}}$. Furthermore, the $B$-quadreal of 
$\overline{p_{1}p_{2}p_{3}}$, which is%
\begin{equation*}
l_{B}=\allowbreak \frac{169}{394}\times \left( -\frac{4901}{47\,083}\right)
\times \frac{197}{116}=-\frac{28\,561}{376\,664}
\end{equation*}%
is the $B$-projective quadrea of $\overline{r_{1}r_{2}r_{3}},$ and the $B$%
-quadreal of $\overline{r_{1}r_{2}r_{3}}$ is the $B$-projective quadrea of $%
\overline{p_{1}p_{2}p_{3}}$.

\section{Using general metrics for linear algebra problems}

We make just a few simple observations that give more motivation for the
utility of having a general trigonometry valid for an arbitrary quadratic
form. Suppose we are working in three-dimensional Euclidean space $\mathbb{E}%
^{3}$, and we have an application that crucially involves a linear change of
coordinates, given by a linear transformation $\mathbb{E}:\mathbb{V}%
^{3}\rightarrow \mathbb{V}^{3}$ represented by a matrix of the same name.

In classical geometry, we may find it awkward to make the transition to this
new change of coordinates if we are primarily interested in metrical
properties, for the simple reason that 
\begin{equation*}
\left( v_{1}L\right) \cdot \left( v_{2}L\right) =\left( v_{1}L\right) \left(
v_{2}L\right) ^{T}=v_{1}Bv_{2}^{T}=v_{1}\cdot _{B}v_{2}
\end{equation*}%
involves a new scalar product associated to the matrix 
\begin{equation*}
B\equiv LL^{T}.
\end{equation*}

But with the set up of this paper, we are completely in control of such a
general metrical situation, so we may apply the desired linear
transformation without weakening our ability to apply the powerful
geometrical tools provided by the $B$-vector calculus.

In effect we are introducing the metric as a variable quantity into our
geometrical three-dimensional theories. This is a familiar approach in
modern differential geometry, but it has seen little development in
classical geometry.

As a simple example, suppose we are interested in the geometry of a lattice
in three dimensions and associated theta functions. It is traditional to
consider different lattices, but always with respect to a Euclidean
quadratic form. With this more general technology, we may perform a linear
transformation so that the lattice itself is just the standard lattice, and
all the complexity and variability is then inherent in the now variable
quadratic form.

\section{Further applications to the trigonometry of a tetrahedron}

The affine and projective triangle geometries are both necessary ingredients
for the systematic rational investigation of a tetrahedron in
three-dimensional affine/vector space. As an elementary example, we may use
the $B$-quadrance to develop a new metrical quantity associated to a general
tetrahedron itself. We can also use the projective tools developed in this
paper to develop two new affine metrical quantities: a $B$-dihedral spread
and a $B$-solid spread. These two quantities extend the projective notions
of $B$-projective spread and $B$-projective quadrea respectively into the
three-dimensional affine space. We will deal with these quantities and their
results in a more formal way in the subsequent paper \textit{Generalised
vector products applied to affine rational trigonometry of a general
tetrahedron}.


\begin{thebibliography}{99}
\bibitem{Binet} Binet, J. F. M.: M\'{e}moire sur un systeme de formules
analytiques, et leur application \`{a} des considerations g\'{e}om\'{e}%
triques. \textit{Journal \'{E}cole Polytechnique}, 9, 280-302 (1812).

\bibitem{BS} Brualdi, R. A. \& Schneider, H.: Determinantal identities:
Gauss, Schur, Cauchy, Sylvester, Kronecker, Jacobi, Binet, Laplace, Muir,
and Cayley. \textit{Linear Algebra Appl.}, 52-53, 769-791 (1983).

\bibitem{Cauchy} Cauchy, A.: Memoire sur les fonctions qui ne peuvent
obtenir que deux valeurs \'{e}gales et des signes contraires par suite des
transpositions op\'{e}r\'{e}es entre les variables qu'elles renferment. 
\textit{Journal \'{E}cole Polytechnique}, 17, 29-112 (1815).

\bibitem{CM} Chapman, S. \& Milne, E. A.: The proof of the formula for the
vector triple product. \textit{Math. Gaz.}, 23(253), 35-38 (1939).

\bibitem{Collomb} Collomb, C.: \textit{A tutorial on inverting 3\ by 3\
matrices with cross products}. http://www.emptyloop.com/technotes (n.d.).
Accessed 29 June 2019.

\bibitem{Doran et al} Doran, C. J. L. \& Lasenby, A. N.: \textit{Geometric
Algebra for Physicists}. Cambridge University Press, Cambridge, UK (2003).

\bibitem{Gibbs} Gibbs, J. W. \& Wilson, E. B.: \textit{Vector Analysis: A
text-book for the use of students of mathematics and physics founded upon
the lectures of J. Willard Gibbs}. Yale University Press, New Haven, CT
(1901).

\bibitem{Heath} Heath, T. L.: \textit{The Thirteen Books of Euclid's Elements%
}. Dover Publishing Inc, New York, NY (1956).

\bibitem{Heaviside} Heaviside, O.: \textit{Electromagnetic Theory} (Vol 1).
"The Electrician" Printing and Publishing Company Ltd, London, UK (1894).

\bibitem{Hestenes et al} Hestenes, D. \& Sobczyk, G.: \textit{Clifford
Algebra to Geometric Calculus, a Unified Language for Mathematics and Physics%
}. D. Reidel Publishing Company, Dordrecht, Netherlands (1984).

\bibitem{Jacobi} Jacobi, C. G. J.: \textit{Fundamenta Nova Theoriae
Functionum Ellipticarum}. Pontrieu \& Co. and Treuttel \& Wuerz, Paris,
France (1829).

\bibitem{Lagrange} Lagrange, J. L.: Solutions analytiques de quelques probl%
\`{e}mes sur les pyramides triangulaires. \textit{Oeuvres de Lagrange}, 
\textbf{3}, 661-692 (1773).

\bibitem{Minkowski} Minkowski, H.: Das Relativit\"{a}tsprinzip. \textit{Ann.
Phys.}, 352(15), 927-938 (1907).

\bibitem{Narayan} Narayan, S.: \textit{A Textbook of Vector Analysis}. S.
Chand \& Company Ltd, New Delhi, India (1961).

\bibitem{Spiegel} Spiegel, M. R.: \textit{Schaum's Outline of Vector Analysis%
}. McGraw-Hill Inc, New York, NY (1959).

\bibitem{WildAP} Wildberger, N. J.: \textit{Affine and projective universal
geometry}. https://arxiv.org/abs/math/0612499 (2006). Accessed 20 March 2019.

\bibitem{WildDP} Wildberger, N. J.: \textit{Divine Proportions: Rational
Trigonometry to Universal Geometry}. Wild Egg Books, Sydney, Australia
(2005).

\bibitem{WildUHG1} Wildberger, N. J.: Universal hyperbolic geometry I:
trigonometry. \textit{Geom. Dedicata}, 163(1), 215-274 (2013).
\end{thebibliography}
\end{document}